\documentclass[11pt,reqno,twoside]{amsart}

\usepackage{graphicx}
\usepackage{amssymb}
\usepackage{epstopdf}
\usepackage{geometry}
\usepackage{booktabs} % for much better looking tables
\usepackage{array} % for better arrays (eg matrices) in maths
\usepackage{paralist} % very flexible & customisable lists (eg. enumerate/itemize, etc.)
\usepackage{verbatim} % adds environment for commenting out blocks of text & for better verbatim
\usepackage{subfig} % make it possible to include more than one captioned figure/table in a single float
\usepackage{tabularx}
\usepackage{amsmath,amsfonts,amsthm,mathrsfs,amssymb,cite}
\usepackage[usenames]{color}
\usepackage{bm}

\newtheorem{thm}{Theorem}[section]

\newtheorem{lem}{Lemma}[section]
\newtheorem{prop}{Proposition}[section]
\theoremstyle{definition}
\newtheorem{defn}{Definition}[section]
\theoremstyle{remark}

\newtheorem{rem}{Remark}[section]
\numberwithin{equation}{section}

\newcommand{\bu}{\mathbf{u}}

\newcommand{\bvarphi}{\bm{\varphi}}
\newcommand{\bpsi}{\bm{\psi}}

\newcommand{\bnu}{\bm{\nu}}
\newcommand{\bGa}{\mathbf{\Gamma}}
\newcommand{\bx}{\mathbf{x}}
\newcommand{\by}{\mathbf{y}}
\newcommand{\rmi}{\mathrm{i}}
\newcommand{\bS}{\mathbf{S}}
\newcommand{\bK}{\mathbf{K}}
\newcommand{\bI}{\mathbf{I}}
\newcommand{\bbS}{\mathbb{S}}

\newcommand{\bc}{\mathbf{c}}

\newcommand{\bp}{\mathbf{p}}

\newcommand{\ri}{\mathrm{i}}
\newcommand{\bxi}{\bm{\xi}}
\newcommand{\hbxi}{\hat{\bm{\xi}}}
\newcommand{\bzeta}{\bm{\zeta}}
\newcommand{\hbze}{\hat{\bm{\zeta}}}
\newcommand{\bal}{\bm{\alpha}}
\newcommand{\bbe}{\bm{\beta}}

\newcommand{\norm}[1]{\left\Vert#1\right\Vert}
\newcommand{\td}{\tilde}

\newcommand{\Ical}{\mathcal{I}}
\newcommand{\Acal}{\mathcal{A}}
\newcommand{\Bcal}{\mathcal{B}}

\newcommand{\Lcal}{\mathcal{L}}

\newcommand{\Gcal}{\mathcal{G}}
\newcommand{\Hcal}{\mathcal{H}}
\newcommand{\Mcal}{\mathcal{M}}
\newcommand{\Ncal}{\mathcal{N}}
\newcommand{\Ocal}{\mathcal{O}}
\newcommand{\Tcal}{\mathcal{T}}

\newcommand{\Pcal}{\mathcal{P}}

\usepackage[hidelinks]{hyperref}
\title[Dipolar resonances of hard inclusions within a soft elastic material]{Mathematical theory on dipolar resonances of hard inclusions within a soft elastic material}

\author{Hongjie Li}
\address{Department of Mathematics, The Chinese University of Hong Kong, Shatin, Hong Kong SAR, China.}
\email{hjli@math.cuhk.edu.hk}

\author{Jun Zou}
\address{Department of Mathematics, The Chinese University of Hong Kong, Shatin, Hong Kong SAR, China.
}
\email{zou@math.cuhk.edu.hk}

\begin{document}

\maketitle

\begin{abstract}

The current study is motivated by the paper [Z. Liu, et al., {\it Science}, 289(5485), 2000], which investigates the incorporation of hard inclusions within a soft elastic matrix (HISE). The objective is to attain a negative mass density, which is caused by sub-wavelength dipolar resonances. This paper offers a comprehensive and mathematically rigorous understanding of the dipolar resonances within the HISE structure. Firstly, the original formula for the resonant frequencies of arbitrarily shaped hard inclusions is derived explicitly for the first time. Additionally, we demonstrate that the resonances arise due to the high contrast in Lamé parameters between various materials, while the sub-wavelength characteristics of the resonances are attributed to the significant disparities in material densities.
 Furthermore, the dipolar characterization of the resonances is rigorously justified. 
 To validate our findings, we investigate resonant phenomena within a spherical geometry, offering both an alternative method for calculating resonant frequencies and a confirmation of the results presented in this paper.

\medskip



\end{abstract}

\section{Introduction}

Negative-Index material (NIM) is a kind of metamaterials. The NIMs could exhibit a negative index of refraction within a certain frequency range,  a characteristic not found in natural materials. Due to their exotic properties, the metamaterials have attracted significant attention over the past few decades and possess many important applications, including invisibility cloaking \cite{ACK6055, LL0165, LL9023, LLL6014, LL6090, MN1715}, super-resolution imaging \cite{CLH1677, BLLW9091, LLLW9004} and waveguides \cite{AE1274, AHY4171, LP5021}. However, it is essential to note that such metamaterials do not occur naturally and must be artificially engineered.  The creation of NIMs typically involves the deliberate design of sub-wavelength resonators, often referred to as ``meta-atoms''. The term ``sub-wavelength'' implies that the dimensions of these meta-atoms are significantly smaller than the wavelength of the incident waves, while ``resonators'' indicate that these meta-atoms resonate with incident waves. In the realm of acoustic metamaterials, the meta-atoms are often realized as bubbles embedded in liquids \cite{AFGL007, CP7599, PhReE.7005}. For electromagnetic metamaterials, silicon nanoparticles, split-ring resonators, and conducting wires are commonly employed as meta-atoms \cite{ALLZ3224, SC4846, KML2472}.

However, it is worth noting that elastic metamaterials are less touched compared to their counterparts in acoustic and electromagnetic metamaterials. This disparity is primarily attributable to the inherent complexity of elastic materials, where longitudinal and transverse waves are coupled, giving rise to a diverse range of behaviors. As a result, elastic waves exhibit a unique blend of optical and acoustic characteristics. In particular, the elastic material is characterized by three essential physical parameters, i.e. the compression modulus, mass density and shear modulus. To attain negative values for these parameters, metamaterials must induce specific resonances, including monopolar, dipolar, and quadrupolar resonances, which correspond to the compression modulus, mass density, and shear modulus, respectively \cite{NaMLai2011}.
Previous investigations, such as in \cite{CTL1646, LSM0301}, have explored the use of bubbles embedded in soft elastic materials (referred to as BSE structures) to generate monopolar resonance and achieve a negative compression modulus.
 Notably, the paper \cite{LLZ0572} provided the initial mathematical proof for the feasibility of BSE structures.
To obtain the negative shear modulus, the quadrupolar resonances are induced by different structures, including the rubber-coated epoxy spheres and steel-coated soft silicone rubber \cite{PhRevB7919, NaMLai2011}. 
The primary focus of the current study is the utilization of hard inclusions embedded in soft elastic materials (referred to as HISE structures) to induce the dipolar resonance, which therefore could induce the negative mass density \cite{Liu00scie, PhReB7101}. 

The contributions of the paper are as follows. Firstly, we explicitly derive the original formula for the resonant frequencies of arbitrarily shaped hard inclusions for the HISE structure within Theorem \ref{thm:ref}. To the best of our knowledge, this is the first result on the explicit expressions of the resonant frequencies for the HISE structure. Moreover, we characterize the resonance intensively and deeply. We theoretically reveal that the resonance is caused by the high contrast of the Lam\'e parameters (the compression modulus and the shear modulus)  within the hard inclusions and the surrounding elastic matrix, while the sub-wavelength characteristic of the resonance is determined by the contrast of the densities. Indeed, in most of the physical references, only the  high contrast of the Lam\'e parameters is emphasized to achieve the resonance, e.g. \cite{PhReB7101}. Thus the resonant frequencies obtained are not so deep since the contrast of the densities is ignored. Consequently, to have sub-wavelength resonances, both high contrasts of the Lam\'e parameters and the densities should be taken into consideration, which provides insightful instructions to design the NIMs in elasticity. Please refer to Remark \ref{rem:depre} for more details. Secondly, the dipolar characterization  of the resonances is rigorously justified in Theorem \ref{thm:dir}. When the frequencies of the incident wave are located in different regimes, the wave fields inside the hard inclusion are derived explicitly. 
Thirdly, the resonant frequencies for the HISE structure are derived exactly for the spherical geometry using two different methods. One is to utilize the conclusion obtained in Theorem \ref{thm:ref}. The other is through directly solving the original problem to obtain the resonant frequencies. These results are in accordance with each other, thus validating the conclusions drawn in Theorem \ref{thm:ref}. At last, we would like to mention that the results obtained in the current study can be extended to the two dimensional case following a similar discussion. 

The paper is organized as follows. In Section 2, we present the mathematical formulation of our study. Section 3 is devoted to the general exposition of the configuration, and some preliminaries for our future analysis. In Section 4, we derive the explicit expressions of the resonant frequencies, and substantiate the dipolar characterization of the resonances. Section 5 is devoted to the dipolar resonances for the spherical geometry.

\section{Mathematical formulation}

In this section, we present the mathematical formulation for our study.
Let $D$ be a bounded domain in $\mathbb{R}^3$ with connected Lipschitz boundaries. 
Suppose that the background $\mathbb{R}^3\backslash\overline{D}$ is occupied by a regular elastic material parameterized by Lam\'e constants $(\lambda,\mu)$ satisfying the following strong convexity condition 
 \begin{equation}\label{eq:con}
  \mathrm{i)}.~~\mu>0\qquad\mbox{and}\qquad \mathrm{ii)}.~~3\lambda+2\mu>0.
 \end{equation}
 Let $\rho\in\mathbb{R}_+$ denote the density in the background $\mathbb{R}^3\backslash\overline{D}$.
%Consider the case that the domain $D$ is a hard inclusion. 
The corresponding Lem\'e parameters in the domain $D$ are $(\td{\lambda}, \td{\mu})$ satisfying the condition \eqref{eq:con} and the density is $\td{\rho}\in\mathbb{R}_+$.
Let $\bu^i$ be a time-harmonic incident elastic wave satisfying the elastic equation in the entire space $\mathbb{R}^3$
\begin{equation}\label{eq:inci}
\mathcal{L}_{ {\lambda}, {\mu}}\bu^i(\bx) + \omega^2  {\rho} \bu^i(\bx) =0,
\end{equation}
where $\omega>0$ is the angular  frequency. 
In \eqref{eq:inci}, the Lam\'e operator $ \mathcal{L}_{\lambda, \mu}$ associated with the parameters $(\lambda,\mu)$ is defined by 
\begin{equation}\label{op:lame}
 \Lcal_{\lambda,\mu}\bu^i:=\mu \triangle\bu^i + (\lambda+ \mu)\nabla\nabla\cdot\bu^i.
\end{equation}
Then the total displacement field $\bu$ described above is controlled by the following system
\begin{equation}\label{eq:mo}
  \left\{
    \begin{array}{ll}
        \mathcal{L}_{\td{\lambda}, \td{\mu}}\bu(\bx) + \omega^2 \td{\rho} \bu(\bx) =0  &  \bx\in D,  \medskip\\
       \mathcal{L}_{\lambda, \mu}\bu(\bx) + \omega^2\rho\bu(\bx) =0    &   \bx\in \mathbb{R}^3\backslash \overline{D},  \medskip \\
      \bu(\bx)|_- = \bu(\bx)|_+      & \bx\in\partial D,  \medskip \\
      \partial_{\td{\bnu}}\bu(\bx)|_- = \partial_{{\bnu}}\bu(\bx)|_+ & \bx\in\partial D,  \medskip \\
      \bu^s:=\bu-\bu^i  \qquad \mbox{satisfies the radiation condition},
    \end{array}
  \right.
\end{equation} 
where the subscript $\pm$ indicate the limits from outside and inside of $D$, respectively.
In \eqref{eq:mo}, the co-normal derivative $\partial_{\bnu}$ associated with the parameters $(\lambda, \mu)$ is defined by 
\begin{equation}\label{eq:trac}
\partial_{\bnu}\bu=\lambda(\nabla\cdot \bu)\bnu + 2\mu(\nabla^s\bu) \bnu,
\end{equation}
where $\bnu$ represents the outward unit normal to $\partial D$ and the operator $\nabla^s$ is the symmetric gradient
 \begin{equation}\label{eq:sg1}
 \nabla^s\mathbf{u}:=\frac{1}{2}\left(\nabla\mathbf{u}+\nabla\mathbf{u}^t \right),
 \end{equation}
 with $\nabla\bu$ denoting the matrix $(\partial_j u_i)_{i,j=1}^3$ and the superscript $t$ signifying the matrix transpose.
The operators $ \mathcal{L}_{\td{\lambda}, \td{\mu}}$ and $\partial_{\td{\bnu}}$ are defined in \eqref{op:lame} and \eqref{eq:trac} , respectively, with the parameters $(\lambda, \mu)$ replaced by $(\td{\lambda}, \td{\mu})$.
 In \eqref{eq:mo}, the radiation condition designates the following facts \cite{KupTPET, LLL9974}:
\begin{equation}\label{eq:radi}
\begin{split}
(\nabla\times\nabla\times \bu^s)(\bx)\times\frac{\bx}{|\bx|}-\mathrm{i} {k}_s\nabla\times \bu^s(\bx)=&\mathcal{O}(|\bx|^{-2}),\\
\frac{\bx}{|\bx|}\cdot[\nabla(\nabla\cdot \bu^s)](\bx)-\mathrm{i} {k}_p\nabla \bu^s(\bx)=&\mathcal{O}(|\bx|^{-2}),
\end{split}
\end{equation}
as $|\mathbf{x}|\rightarrow+\infty$, where $\rmi$ signifies the imaginary unit and
\begin{equation}\label{pa:ksp}
 {k}_s=\omega/{c}_s, \quad  \quad {k}_p=\omega/{c}_p,
\end{equation}
with
\[
{c}_s = \sqrt{{\mu}/{\rho}}, \quad \quad {c}_p=\sqrt{ ({\lambda} + 2 {\mu})/{\rho}}.
\]
It is well-known that the elastic wave can be decomposed into the shear wave (s-wave) and the compressional wave (p-wave). In \eqref{pa:ksp}, the parameters $k_s$ and $k_p$ signify the wavenumbers of the s-wave and p-wave, respectively. 
In what follows, the parameters $\td{k}_s, \td{k}_p, \td{c}_s, \td{c}_p$ are defined in \eqref{pa:ksp} by replacing $(\lambda, \mu, \rho)$ with $(\td{\lambda}, \td{\mu}, \td{\rho})$.

In this paper, we mainly apply the potential theory to investigate the system \eqref{eq:mo}. To that end, we first introduce the potential theory of the Lam\'e system.
The fundamental solution $\bGa^{\omega}=(\Gamma^{\omega}_{i,j})_{i,j=1}^3$ of the operator $\Lcal_{\lambda,\mu} + \rho\omega^2$ in three dimensions is given by \cite{DLL7678}
\begin{equation}\label{eq:ef}
 (\Gamma^{\omega}_{i,j})_{i,j=1}^3(\bx)=-\frac{\delta_{ij}}{4\pi\mu|\bx|}e^{\rmi k_s |\bx|} + \frac{1}{4\pi \omega^2\rho}\partial_i\partial_j\frac{e^{\rmi k_p|\bx|} - e^{\rmi k_s|\bx|}}{|\bx|},
\end{equation}
where $\delta_{ij}$ is Kronecker delta function. In particular, when $\omega=0$, we denote  $\bGa^{0}$ by  $\bGa$ for simplicity, and $\bGa$ has the following expression
\[
 \Gamma_{i,j} (\bx)= -\frac{1}{8\pi}\left( \frac{1}{\mu} + \frac{1}{\lambda + 2\mu} \right)  \frac{\delta_{ij}}{|\bx|}  -\frac{1}{8\pi}\left( \frac{1}{\mu} - \frac{1}{\lambda + 2\mu} \right)  \frac{\bx_i \bx_j}{|\bx|^3 }. 
\]
Then the single-layer potential associated with the fundamental solution $\bGa^{\omega}$ is defined by
\begin{equation}\label{eq:single}
 \bS_{\partial D}^{\omega}[\bvarphi](\bx)=\int_{\partial D} \bGa^{\omega}(\bx-\by)\bvarphi(\by)ds(\by), \quad \bx\in\mathbb{R}^3,
\end{equation}
for $\bvarphi\in L^2(\partial D)^3$. On the boundary $\partial D$, the conormal derivative of the single-layer potential satisfies the following jump formula
\begin{equation}\label{eq:jump}
\partial_{\bnu}   \bS_{\partial D}^{\omega}[\bvarphi]|_{\pm}(\bx)=\left( \pm\frac{1}{2}\bI +  \bK_{\partial D}^{\omega, *} \right)[\bvarphi](\bx) \quad \bx\in\partial D,
\end{equation}
where
\[
  \bK_{\partial D}^{\omega, *} [\bvarphi](\bx)=\mbox{p.v.} \int_{\partial D} \partial_{\bnu_{\bx}} \bGa^{\omega}(\bx-\by)\bvarphi(\by)ds(\by),
\]
with $\mbox{p.v.}$ standing for the Cauchy principal value. We would like to mention that the operator $ \bK_{\partial D}^{\omega, *}$ in \eqref{eq:jump} is called the Neumann-Poincar\'e (N-P) operator, which is a critical operator in the investigation of metamaterials in elasticity. In what follows, we denote $\bS_{\partial D}^{0}$, $\bK_{\partial D}^{0, *}$ by $\bS_{\partial D}$, $\bK_{\partial D}^{ *}$, respectively, for simplicity. 

With the help of the potential theory presented above, the solution to the system \eqref{eq:mo} can be written as 
\begin{equation}\label{eq:sol}
  \bu=
 \left\{
   \begin{array}{ll}
     \td{\bS}_{\partial D}^{\omega}[\bvarphi](\bx), & \bx\in D,  \smallskip \\
     {\bS}_{\partial D}^{\omega}[\bpsi](\bx) +\bu^i, &  \bx\in \mathbb{R}^3\backslash \overline{D},
   \end{array}
 \right.
\end{equation}
for some density functions $\bvarphi, \bpsi \in L^2(\partial D)^3$. Here, the operator $\td{\bS}_{\partial D}^{\omega}$ is the singley-layer potential operator associated with the parameters $(\td{\lambda}, \td{\mu})$ and $\td{\rho}$. 
By matching the transmission conditions on the boundary, i.e. the third and fourth conditions in \eqref{eq:mo} and with the help of the jump formula in \eqref{eq:jump}, the density functions $\bvarphi$ and $\bpsi$ in \eqref{eq:sol} should satisfy the following system:
\begin{equation}\label{eq:or}
\Acal(\omega,\delta) [\Phi](\bx)=F(\bx), \quad \bx\in\partial D,
\end{equation}
where
\[
\Acal(\omega,\delta)=  \left(
    \begin{array}{cc}
       \td{\bS}_{\partial D}^{\omega} &  -{\bS}_{\partial D}^{\omega}\medskip \\
      -\frac{I}{2} +  \td{\bK}_{\partial D}^{\omega, *} & -\frac{I}{2} -  {\bK}_{\partial D}^{\omega, *}\\
    \end{array}
  \right),
  \;\;
  \Phi= \left(
    \begin{array}{c}
      \bvarphi \\
     \bpsi \\
    \end{array}
  \right)
   \;\; \mbox{and} \;\; 
  F= \left(
    \begin{array}{c}
     \bu^i \\
     \partial_{\bnu}  \bu^i \\
    \end{array}
  \right).
\]
For the further discussion, we introduce the spaces $\Hcal=L^2(\partial D)^3\times L^2(\partial D)^3$ and $\Hcal^1=H^1(\partial D)^3\times L^2(\partial D)^3$. Apparently, the operator $\Acal(\omega,\delta)$ is defined from $\Hcal$ to $\Hcal^1$. Next, we give the definition of the sub-wavelength resonance of the scattering system \eqref{eq:mo} based on the operator $\Acal(\omega,\delta)$.
\begin{defn}
The sub-wavelength resonance of the scattering system \eqref{eq:mo} occurs if there exists a frequency $\omega\ll 1$ such that the operator $\Acal(\omega,\delta)$ has a nontrivial kernel, i.e. 
\begin{equation}\label{eq:conre1}
 \Acal(\omega,\delta)[\Phi](\bx)=0,
\end{equation}
for some nontrivial $\Phi\in\Hcal$.
\end{defn}

%
%
%The resonance of the scattering problem \eqref{eq:mo} can be defined for some frequency $\omega$ such that there exists a nontrivial solution $\Phi$ to the following equation 
%\begin{equation}\label{eq:conre}
% \mfa(\omega,\delta)[\Phi](\bx)=0.
%\end{equation}
%In the rest of the paper, instead of directly considering the equation \eqref{eq:conre}, we consider the following equation 
%\begin{equation}\label{eq:conre1}
% \Acal(\omega,\delta)[\Phi](\bx)=0,
%\end{equation}
%where 
%\[
%\Acal(\omega,\delta)=  \left(
%    \begin{array}{cc}
%       \delta{\bS}_{\partial D}^{\omega} &  -\delta\td{\bS}_{\partial D}^{\omega}\medskip \\
%      -\frac{I}{2} +  \bK_{\partial D}^{\omega, *} & -\frac{I}{2} -  \td{\bK}_{\partial D}^{\omega, *}\\
%    \end{array}
%  \right).
%\]
%
%Indeed, the last two equations share the same solutions as stated in the following lemma.
%\begin{lem}
%The function $\Phi$ satisfies the equation $\mfa(\omega,\delta)[\Phi](\bx)=0$ if and only if the function $\Phi$ satisfies $\Acal(\omega,\delta)[\Phi](\bx)=0$.
%\end{lem}
%\begin{proof}
%The proof is obvious.
%\end{proof}

\section{General requirements of the configuration, preliminaries and auxiliary results}

In this section, we present general requirements of the configuration of the system \eqref{eq:mo}, i.e. the physical parameters in the regions of the domain $D$ and the free space. Then some preliminaries and auxiliary results shall be provided.

In this paper, we mainly consider the configuration that the hard inclusion $D$ is immersed within a soft elastic material, such as lead inclusions coated by the silicone rubber \cite{Liu00scie}. Then the corresponding physical parameters in different regions satisfy the following relationship
\begin{equation}\label{eq:hicon}
    \td{\lambda} = \frac{1}{\delta }\lambda, \quad  \td{\mu} = \frac{1}{\delta}\mu, \quad  \td{\rho} =\frac{1}{\epsilon} \rho,
\end{equation}
where $\delta\ll1$ and $\epsilon\ll 1$. Next, we induce another parameter $\tau$ describing the contrast of the wave speed in different regions, namely 
\begin{equation}\label{eq:detau}
\tau = \frac{{c}_s}{\td{c}_s} = \frac{{c}_p}{\td{c}_p}= \sqrt{\delta/\epsilon}.
\end{equation}
The last equality in \eqref{eq:detau} follows from \eqref{eq:hicon} and the expressions of $c_s$ as well as $c_p$ in \eqref{pa:ksp}. Moreover, we assume that the contrast $\tau$ satisfies 
\begin{equation}\label{eq:de}
  \tau= \sqrt{\delta/\epsilon} \leq \Ocal(1).
\end{equation}
The assumption in \eqref{eq:de} is reasonable since wave speed inside the hard inclusion is larger than that in the soft material. 
%We would like to mention that the assumptions here are reasonable and include almost all the composite materials of hard inclusions embedded into soft elastic materials, such as [...] (physical paper).

In the current study, we consider the subwavelength resonance; that is the size of the domain $D$ is smaller compared with the wavelength of the incident wave. 
By coordinate transformation, we may assume that the size of the domain $D$ is of order $1$ and ${k}_s =o(1), {k}_p = o(1)$.
 From the relationship \eqref{eq:de}, we further have that $\td{k}_s = o(1), \td{k}_p = o(1)$ and satisfy
 \begin{equation}
 \td{k}_s = \tau k_s,  \qquad    \td{k}_p =  \tau k_p.  
 \end{equation}

For the further analysis, we first present the asymptotic expansion of the fundamental solution $\bGa^{\omega}$ defined in \eqref{eq:ef}.
\begin{lem}\label{lem:asmfun}
Suppose $\omega\in\mathbb{R}_+$ and $\omega\ll 1$, then there holds the following asymptotic expansion of the fundamental solution $\bGa^{\omega}$ defined in \eqref{eq:ef}
\begin{equation}\label{eq:fse}
  \mathbf{\Gamma}^{\omega}(\bx) = \sum_{n=0}^{\infty} \omega^n  \mathbf{\Gamma}_n(\bx),
\end{equation}
where
\[
 \begin{split}
  \mathbf{\Gamma}_n(\bx) = & -\frac{\alpha_1}{4\pi} \frac{\mathrm{i}^n}{(n+2)n!}\left(\frac{n+1}{c^{n}_s} + \frac{1}{c^{n}_p} \right)  |\bx|^{n-1}\mathbf{I} \\
     & +\frac{\alpha_2}{4\pi}\frac{\mathrm{i}^n(n-1)}{(n+2)n!}\left(\frac{1}{c^{n}_s} - \frac{1}{c^{n}_p} \right)  |\bx|^{n-3}\bx\bx^T,
 \end{split}
\]
with 
\[
\alpha_1= \left( \frac{1}{\mu} + \frac{1}{2\mu +\lambda} \right) \quad \mbox{and} \quad \alpha_2 = \left( \frac{1}{\mu} - \frac{1}{2\mu +\lambda} \right).
\]
\end{lem}
\begin{proof}
For $k\ll 1$, the following function enjoys the asymptotic expansion
\begin{equation*}
 \displaystyle{-\frac{e^{\mathrm{i} k|\bx|}}{4\pi |\bx|}}=-\sum_{j=1}^{\infty} \frac{\rmi}{4\pi} \frac{k^j(\rmi|\bx|)^{j-1}}{j!}.
\end{equation*}
Substituting the last expansion into the expression of the fundamental solution \eqref{eq:ef} and by tidies calculation, one can obtain \eqref{eq:fse}.
The proof is completed.
\end{proof}

From Lemma \ref{lem:asmfun}, the single layer potential $\bS_{\partial D}^{\omega}$ defined in \eqref{eq:single} has the following asymptotic expansion:
\begin{equation}\label{eq:sise}
\bS_{\partial D}^{\omega} = \bS_{\partial D} + \sum_{n=1}^\infty \omega^n  \bS_{\partial D, n},
\end{equation}
where
\[
 \bS_{\partial D, n}[\bvarphi](\bx)=\int_{\partial D} \bGa_n(\bx-\by)\bvarphi(\by)ds(\by).
\]
In particular, one has that
\begin{equation}\label{eq:s1}
\bS_{\partial D, 1}[\bvarphi](\bx)= -\frac{\ri \alpha_1}{12\pi} \left(\frac{2}{c_s} + \frac{1}{c_p} \right) \int_{\partial D} \bvarphi(\by)ds(\by).
\end{equation}

\begin{lem}
The norm $\norm{\bS_{\partial D, n}}_{\Lcal(L^2(\partial D), H^1(\partial D))} $ is uniformly bounded with respect to $n$. Moreover, the series in \eqref{eq:sise} is convergent in $\Lcal(L^2(\partial D), H^1(\partial D))$ (c.f. \cite{DLL1067}).
\end{lem}

Next, we consider the asymptotic expansion for the Neumann-Poincar\'e (N-P) operator $ \bK_{\partial D}^{\omega, *}$. Using Lemma \ref{lem:asmfun} yields the following asymptotic expansion for the N-P operator $ \bK_{\partial D}^{\omega, *}$:
\begin{equation}\label{eq:npse}
\bK_{\partial D}^{\omega, *} = \bK_{\partial D}^{ *} + \sum_{n=1}^\infty \omega^n  \bK_{\partial D, n}^{ *},
\end{equation}
where
\[
 \bK_{\partial D, n}^{ *} [\bvarphi](\bx)=\int_{\partial D}  \frac{\partial \bGa_n}{\partial \bnu(\bx)} (\bx-\by)\bvarphi(\by)ds(\by).
\]
In particular, we have 
\[
 \bK_{\partial D, 1}^{ *}=0,
\]
which follows from that the function $ \mathbf{\Gamma}_1(\bx)$ is a constant. Moreover, the follow lemma holds.

\begin{lem}
The norm $\norm{  \bK_{\partial D, n}^{ *}  }_{\Lcal(L^2(\partial D))} $ is uniformly bounded with respect to $n$. Moreover, the series in \eqref{eq:npse} is convergent in $\Lcal(L^2(\partial D))$ (c.f. \cite{DLL1067}).
\end{lem}

Then we introduce the vector space $\mho$ that is spanned by all linear solutions to the equation
\begin{equation}
  \left\{
    \begin{array}{ll}
        \mathcal{L}_{ {\lambda}, {\mu}}\bu=0  &  \bx\in D,  \medskip\\
      \partial_{{\bnu}}\bu=0     & \bx\in\partial D.  \medskip
    \end{array}
  \right.
\end{equation} 
 Indeed, the space $\mho$ can be explicitly expressed as \cite{ABJH6625}
\begin{equation}
\begin{split}
\mho = &\left\{  \mathbf{a} + \mathbf{B}\bx, \mathbf{a}\in\mathbb{R}^3, \mathbf{B}\in M^A \right\}, \\
%= & \mbox{span} \left\{  \bxi_1,  \bxi_2,  \bxi_3,  \bxi_4,  \bxi_5,  \bxi_6 \right\} ,
\end{split}
\end{equation}
where $M^A$ is the space of antisymmetric matrices with the size $3\times 3$. Direct calculation shows that the dimension of the space $\mho$ is $6$ and the space $\mho$ is spanned by 
\begin{equation}\label{eq:dpsi}
\left[\begin{array}{l}
1 \\
0\\
0
\end{array}\right],\;\; \left[\begin{array}{l}
0 \\
1 \\
0
\end{array}\right],\;\; \left[\begin{array}{c}
0 \\
0 \\
1
\end{array}\right], \;\;\left[\begin{array}{c}
x_2 \\
-x_1 \\
0
\end{array}\right],\;\; \left[\begin{array}{c}
x_3 \\
0 \\
-x_1
\end{array}\right],\;\; \left[\begin{array}{c}
0 \\
x_3 \\
-x_2
\end{array}\right].
\end{equation}

\begin{lem}\label{lem:kek}
The kernel of the operator $-1/2 +  \bK_{\partial D}$ coincides with the space $\mho$, where $\bK_{\partial D}$ is the adjoint operator of $ \bK^*_{\partial D}$ (c.f. \cite{ABJH6625}).
\end{lem}

 Denote by $\mho^*$ the kernel of the operator $-1/2 +  \bK^*_{\partial D}$. The following facts can be found in \cite{DKV7353}. 
\begin{lem}\label{lem:kekske}
The dimension of the kernel of the operator $-1/2 +  \bK^*_{\partial D}$ is $6$, i.e. $\dim(\mho^*)=6$.
\end{lem}

Then we define the space $L^2_{\bxi}(\partial D)^3$ by 
\begin{equation}
L^2_{\bxi}(\partial D)^3 :=\left\{ \bvarphi\in L^{2}(\partial D)^3:  \int_{\partial D} \bvarphi \cdot \bxi =0, \quad \forall\; \bxi\in\mho  \right\}.
\end{equation}
One has that 
\begin{lem}\label{lem:keks}
The operator $-1/2 +  \bK^*_{\partial D}$ is invertible in the space $L^2_{\bxi}(\partial D)^3$ (c.f. \cite{DKV7353}).
\end{lem}

The following lemma shall be used in our future analysis.
\begin{lem}\label{lem:orker}
There exist basis functions $\{\bxi_j\}^6_{j=1}$ of $\mho$ and $\{\bzeta_j\}^6_{j=1}$ of $\mho^*$ such that the following properties hold
\begin{enumerate}[(i)]
\item $(\bzeta_i, \bzeta_j)=\delta_{ij}$,
\item $(\bzeta_i, \bxi_j)=\delta_{ij}$,
\item $(\bxi_i, \bxi_j)=c_i \delta_{ij}$,
\end{enumerate}
 with $c_i\in\mathbb{R}_+$, where $(\cdot, \cdot)$ denotes the $L^2(\partial D)$ inner product. 
\end{lem}

\begin{proof}
We first choose $\{\hat{\bxi}_j\}^6_{j=1} \subset \mho$ such that $(\hbxi_i, \hbxi_j)=\delta_{ij}$. Direct calculation shows that $ \partial_{\bnu} \bS_{\partial D}[\hbxi_i ] |_{-} \in L^2_{\bxi}(\partial D)^3$. Lemma \ref{lem:keks} shows that the operator $1/2 -  \bK^*_{\partial D}$ is invertible on $L^2_{\bxi}(\partial D)^3$. Thus there exists a unique $\hbze_i$ satisfying 
 \[
  \left( 1/2 -  \bK^*_{\partial D} \right) [\hbze_i] =  \partial_{\bnu} \bS_{\partial D}[\hbxi_i ] |_{-} =  \left( -1/2 +  \bK^*_{\partial D} \right) [\hbxi_i ]. 
 \]
Define $\tilde\bzeta_i = \hbze_i + \hbxi_i$ and then the last equation gives that 
\begin{equation}
 \left( -1/2 +  \bK^*_{\partial D} \right)[\tilde\bzeta_i]=0. 
\end{equation}
Moreover, the following relationship holds
\[
 (\tilde\bzeta_i, \hat\bxi_j) =  ( \hbze_i + \hbxi_i, \hat\bxi_j) =  ( \hbze_i, \hat\bxi_j) +( \hbxi_i, \hat\bxi_j) = \delta_{ij}.
\]
If the vectors $\tilde\bzeta_i$ satisfy the condition $(\tilde\bzeta_i, \tilde\bzeta_j)=\delta_{ij}$, just let $\bzeta_i=\tilde\bzeta_i$. Otherwise, let the matrix $P$ be defined by $P=(p_{ij})_{i,j=1}^6$, where $p_{ij} = (\tilde\bzeta_i, \tilde\bzeta_j)$. Apparently, the matrix $P$ is real, symmetric and positive-definite and the following eigendecomposition holds
\[
 PT= T \Lambda,
\]
where the matrix $T$ is a orthogonal matrix and $\Lambda$ is a diagonal matrix whose diagonal elements are the eigenvalues of the matrix $P$. Then we induce the transformation matrix $R$ defined by 
\[
 R_{:,i}=\frac{1}{\sqrt{\Lambda_{ii}}} T_{:,i},
\]
where $R_{:,i}$ denotes the i-th column of the matrix $R$. Then the functions defined as follows 
\[
 (\bzeta_1, \bzeta_2, \bzeta_3, \bzeta_4, \bzeta_5, \bzeta_6) = (\tilde\bzeta_1, \tilde\bzeta_2, \tilde\bzeta_3, \tilde\bzeta_4, \tilde\bzeta_5, \tilde\bzeta_6) R,
\]
satisfy the condition $(i)$. 

Furthermore, if we define 
\[
( \bxi_1,  \bxi_2,  \bxi_3,  \bxi_4,  \bxi_5,  \bxi_6)  = (\hat\bxi_1, \hat\bxi_2, \hat\bxi_3, \hat\bxi_4, \hat\bxi_5, \hat\bxi_6) \tilde R,
\]
where the transformation matrix $\tilde R$ is defined by 
\[
  \tilde R_{:,i}= \sqrt{\Lambda_{ii}}  T_{:,i}.
\]
Finally, one could verify that 
\[
 (\bzeta_i, \bxi_j)=\delta_{ij}, \quad (\bxi_i, \bxi_j)=c_i \delta_{ij}.
\]
The proof is completed.
\end{proof}

\begin{lem}\label{lem:ma1}
The matrix $M=(m_{ij})_{i,j=1}^6$ is invertible, where
\begin{equation}\label{eq:mij}
m_{ij}=\int_{D} -\bxi_j(\by) \cdot  \bS_{\partial D}[\bzeta_i](\by) d\by,
\end{equation}
where the functions $\bxi_j, \bzeta_j$ with $1\leq j\leq 6$ are given in Lemma \ref{lem:orker}.
\end{lem}
\begin{proof}
    First note that the single-layer potential operator $\bS_{\partial D}$ is invertible from $\mho^*$ to $\mho$ \cite{DKV7353}. Suppose that the matrix $M$ is not invertible. Then there exists a nontrivial vector $\mathbf{a}=\{a_i\}_{i=1}^6$ such that $M\mathbf{a}=0$. Thus there exists a nontrivial function $\tilde{\bxi}=\sum_{i=1}^6 a_i  \bS_{\partial D}[\bzeta_j]$ such that 
    \[
    \int_{D} \bxi_i \cdot  \tilde{\bxi} d\by=0, \quad 1\leq i\leq 6.
    \]
    Then one can conclude that $\tilde{\bxi}=0$, which further shows that $\mathbf{a}=0$ thanks to the invertibility of the operator $\bS_{\partial D}$. The proof is completed.
\end{proof}

At the end of this section, we introduce the generalized Rouch\'e's theorem, which mainly follows from the monograph \cite{AKL7848}.

Let $\mathcal C$ and $\mathcal C^{\prime}$ denote two Banach spaces, and $\mathcal{L}(\mathcal{C}, \mathcal{C}^{\prime})$ be the space of bounded linear operators from $\mathcal C$ into $\mathcal C^{\prime}$.
If $A(z)$ is holomorphic and invertible at $z_0$, then $z_0$ is called a regular point of $A(z)$. The point $z_0$ is called a normal point of $A(z)$ if $A(z)$ is finitely meromorphic, of Fredholm type at $z_0$, and regular in a neighborhood of $z_0$ except at $z_0$ itself.

Let $\mathcal{U}(z_0)$ be the set of all operator-valued functions with values in $\mathcal{L}(\mathcal{C}, \mathcal{C}^{\prime})$  which are holomorphic in some neighborhood of $z_0$, except possibly at $z_0$. The point $z_0$ is called a characteristic value of $A(z)\in\mathcal{U}(z_0)$ if there exists a vector-valued function $\varphi(z)$ with values in $\mathcal{C}$ such that
\begin{enumerate}
\item $\varphi(z)$ is holomorphic at $z_0$ and $\varphi(z_0)\neq 0$;
\item $A(z)\varphi(z)$ is holomorphic at $z_0$ and vanishes at this point.
\end{enumerate}
Here, $\varphi(z)$ is called a root function of $A(z)$ associated with the characteristic value $z_0$. The vector $\varphi_0=\varphi(z_0)$ is called an eigenvector. Then there exists a number $m(\varphi)\geq 1$ and a vector-valued function $\psi(z)$ with values in $\mathcal C^{\prime}$, holomorphic at $z_0$, such that
\[
A(z)\varphi(z)=(z-z_0)^{m(\varphi)}\psi(z), \quad \psi(z_0)\neq 0.
\]
The number $m(\varphi)$ is called the multiplicity of the root function $\varphi(z)$. For $\varphi_0\in\ker A(z_0)$, we define the rank of $\varphi_0$, denoted by $\mbox{rank}(\varphi_0)$, to be the maximum of the multiplicities of all root functions $\varphi(z)$ with $\varphi(z_0)=\varphi_0$.

 Suppose that the space of $\ker A(z_0)$ is spanned by  $\varphi_0^j$, $j=1, \dots, n,$. We call 
 \[
  M(A(z_0)):=\sum_{j=1}^n \mbox{rank}(\varphi_0^j)
 \]
the null multiplicity of the characteristic value $z_0$ of $A(z)$.

Let $V$ be a simply connected bounded domain with rectifiable boundary $\partial V$. An operator-valued function $A(z)$ which is finitely meromorphic and of Fredholm type in $V$ and continuous on $\partial V$ is called normal with respect to $\partial V$ if the operator $A(z)$ is invertible in $V$ , except for a finite number of points of $V$ which are normal points of $A(z)$.

In what follows, we assume that there is no poles for $A(z)$ in $V$. If $A(z)$ is normal with respect to the contour $\partial V$ and $z_i$, $i = 1, . . . , m$, are
all its characteristic values lying in $V$, the full multiplicity $\Mcal (A(z); \partial V) $
of $A(z)$ in $V$ is the number of characteristic values of $A(z)$ in $V$ , counted with
their multiplicities.

\begin{lem}\label{lem:sge}
Let $A(z)$ be an operator-valued function which is normal with respect to $\partial V$. If an operator-valued function $S(z)$, which is finitely meromorphic in $V$ and continuous on $\partial V$ satisfies the condition
\[
\norm{A^{-1}(z) S(z)}_{\Lcal(\mathcal{C}, \mathcal{C})}<1, \quad \forall z\in \partial V,
\]
then $A(z) + S(z)$ is also normal with respect to  $\partial V$ and 
\[
\Mcal (A(z); \partial V) = \Mcal (A(z) + S(z); \partial V),
\]
if there is no poles for $A(z)$ in $V$.
\end{lem}

\section{The derivation of the resonant frequencies and the characterization of the resonances}

In this section, we discuss the sub-wavelength resonant phenomenon  for the system \eqref{eq:mo}. First, the original formula of the resonant frequencies are derived explicitly. Then the dipolar characteristics of the resonances are verified rigorously.
 To that end, we first present the asymptotic expansion for the operator $\Acal(\omega,\delta)$ defined in \eqref{eq:or}.

\begin{lem}\label{lem:deoa}
The operator $\Acal(\omega,\delta)$ defined in \eqref{eq:or} has the following asymptotic expansion
\[
\begin{split}
  \Acal(\omega,\delta)= \Acal_0 + \Bcal(\omega, \delta) = &  \Acal_0 + \omega\Acal_{1,0} + \omega^2\Acal_{2,0} +  \delta\Acal_{0,1}  + \omega^2 \tau^2 \Acal_{0,2} + \\
  & \Ocal \left( \omega^3 + \omega \delta \tau  + \left( \omega \tau \right)^3  \right ), 
\end{split}
\]
where 
\[
\Acal_0 =  \left(
    \begin{array}{cc}
       0 &  -\bS_{\partial D}  \medskip \\
      -\frac{I}{2} +  {\bK}_{\partial D}^{*} & -\frac{I}{2} -  {\bK}_{\partial D}^{ *}\\
    \end{array}
  \right),
 \;\;
 \Acal_{1,0} =  \left(
    \begin{array}{cc}
       0 &  - {\bS}_{\partial D, 1}  \medskip \\
      0 & 0\\
    \end{array}
  \right),
\]
and
\[
\Acal_{2,0} =  \left(
    \begin{array}{cc}
       0 &  -{\bS}_{\partial D,2}  \medskip \\
       0 & -  {\bK}_{\partial D, 2}^{ *}\\
    \end{array}
  \right),
 \;\;
 \Acal_{0,1} =  \left(
    \begin{array}{cc}
       {\bS}_{\partial D}  & 0  \medskip \\
      0 & 0\\
    \end{array}
  \right),
  \;\; 
  \Acal_{0,2} =  \left(
    \begin{array}{cc}
       0 & 0 \medskip \\
       \bK_{\partial D, 2}^{*}  & 0 \\
    \end{array}
  \right).
\]
In the last formula, the parameters $\delta$ and $\tau$ are given in \eqref{eq:hicon} and \eqref{eq:detau}, respectively. 
\end{lem}

\begin{proof}
The proof directly follows from asymptotic expansions of  operators $\bS_{\partial D}^{\omega} $ in \eqref{eq:sise}, $\bK_{\partial D}^{\omega, *}$ in \eqref{eq:npse} and parameters defined in \eqref{eq:hicon} as well as \eqref{eq:de}.
\end{proof}

Let $\Acal_{0}^*$ and ${\bS}_{\partial D}^*$ be the adjoint of $\Acal_{0}$ and ${\bS}_{\partial D}$. One can easily check that the operator $\Acal_{0}^*$ has the following expression
\[
\Acal_0^* =  \left(
    \begin{array}{cc}
       0 &    -\frac{I}{2} +  \bK_{\partial D}   \medskip \\
       -\bS^*_{\partial D} & -\frac{I}{2} -  {\bK}_{\partial D} \\
    \end{array}
  \right),
\]
which is defined from $\Hcal^1$ to $\Hcal$. Then the following lemma holds.
\begin{lem}\label{lem:ker}
We have
\begin{enumerate}
\item  $\ker \Acal_{0}$=span$\{\hat\Phi_i \}$, $1\leq i \leq 6$, where
\begin{equation}\label{eq:ker1}
  \hat\Phi_i= \left(
    \begin{array}{c}
      \bzeta_i \\
     0 \\
    \end{array}
  \right),
\end{equation}
with $\bzeta_i$ given in Lemma \ref{lem:orker}. Moreover, the functions satisfy $\left( \hat\Phi_i,  \hat\Phi_i \right)_{\Hcal}=\delta_{ij}$.

\item $\ker\Acal_{\delta}^*$=span$\{\hat\Psi_i \}$, where
\[
  \hat\Psi_i=d_i \left(
    \begin{array}{c}
      -(\bS^{-1}_{\partial D})^{*}[ \bxi_i] \\
     \bxi_i \\
    \end{array}
  \right),
\]
with $\bxi_i$ given in Lemma \ref{lem:orker} and the constants $d_i$ being chosen such that $\left( \hat\Psi_i,  \hat\Psi_i \right)_{\Hcal^1}=\delta_{ij}$.
\end{enumerate}
\end{lem}

\begin{proof}
The proofs of the two inclusions in the lemma are similar. Here we just provide the proof for the first one and the other can proved following the same discussion. We first calculate the kernel of the operator $\Acal_{0}$, i.e. functions $\hat\Phi$. From the definition of the operator $\Acal_0$ in Lemma \ref{lem:deoa}, one has that the second component of $\hat\Phi$ vanishes due to the invertibility of the operator $\bS_{\partial D}$. Thus, the first component of $\hat\Phi$ should belong to the kernel of the operator $ -I/2 +  {\bK}_{\partial D}^{*}$. Finally, the kernel functions should have expressions in \eqref{eq:ker1}. The statement $\left( \hat\Phi_i,  \hat\Phi_i \right)_{\Hcal}=\delta_{ij}$ follows from the choice of $\bzeta_i$ given in Lemma \ref{lem:orker}. This completes the proof of the first conclusion. 
\end{proof}

\begin{rem}
In Lemma \ref{lem:ker}, the functions $\bxi_i$, $1\leq i\leq 6$, may not be exactly same as these in Lemma \ref{lem:orker}. Nevertheless, since the functions $\bxi_i$, $1\leq i\leq 6$, in Lemma \ref{lem:orker} are linearly independent, one can always find $\{\tilde\bxi_i\}_{i=1}^6$ that is a orthogonal transformation of $\{\bxi_i\}_{i=1}^6$ such that 
\[
d_i \left(  
\left(
    \begin{array}{c}
      -(\bS^{-1}_{\partial D})^{*}[\tilde \bxi_i] \medskip \\
     \tilde\bxi_i \\
    \end{array}
  \right),
  \left(
    \begin{array}{c}
      -(\bS^{-1}_{\partial D})^{*}[\tilde \bxi_i] \medskip \\
     \tilde\bxi_i \\
    \end{array}
  \right)
 \right)_{\Hcal^1}=\delta_{ij},
\]
following the same discussion as that in Lemma \ref{lem:orker}.
To ease the exposition and by abuse of notation, here we still denote that by $\bxi_i$.

\end{rem}

Next, we shall show the existence of the sub-wavelength resonances for the system \eqref{eq:mo}, i.e. there exist nontrivial values $\omega\ll 1$ such that the operator $\Acal(\omega, \delta)$ defined in \eqref{eq:or} has a nontrivial kernel. 
Indeed, the last lemma shows that $\omega=0$ is a characteristic value for the operator-valued analytical function $\Acal(\omega, 0)$. By the Gohberg–Sigal theory presented in the last section, we can have the following existence result for the sub-wavelength resonances of the system \eqref{eq:mo}.

\begin{lem}\label{lem:exres}
Under the assumption $\delta\ll 1$, if the parameter $\epsilon\ll 1$, then there exist characteristic values $\omega^*\ll 1$ for the operator-valued analytic function $\Acal(\omega, \delta)$.
\end{lem}

\begin{proof}
We rewrite the operator $\Acal(\omega,\delta)$ defined in \eqref{eq:or} by 
\[
\Acal(\omega,\delta) = \Gcal_1(\omega,\delta) + \Gcal_2(\omega,\delta),
\]
where 
\[
 \Gcal_1(\omega,\delta)=  \left(
    \begin{array}{cc}
       0 &  -{\bS}_{\partial D}^{\omega}\medskip \\
      -\frac{I}{2} +  \td{\bK}_{\partial D}^{\omega, *} & -\frac{I}{2} -  {\bK}_{\partial D}^{\omega, *}\\
    \end{array}
  \right) \quad
  \mbox{and} \quad
  \Gcal_2(\omega,\delta)= \left(
    \begin{array}{cc}
        \td{\bS}_{\partial D}^{\omega} &  0\medskip \\
     0 & 0\\
    \end{array}
  \right).
\]
Next, we investigate the invertibility of the operator $\Gcal_1(\omega,\delta)$.
Note that the operator $\Gcal_1(\omega,\delta)$ is a Fredholm operator from $\Hcal$ to $\Hcal^1$ and $\omega=0$ is a characteristic value for the operator $\Gcal_1(\omega,\delta)$ from Lemma \ref{lem:ker}. Thus one can find a small curve $\partial V$ enclosing the origin such the operator $\Gcal_1(\omega,\delta)$ is invertible for $\omega\in\partial V$ thanks to the discreteness of the spectrum of the Fredholm operator $\Gcal_1(\omega,\delta)$. Direct calculation shows that for $\omega\in\partial V$
\[
\begin{split}
    \Gcal_1^{-1}(\omega,\delta) = \left(
    \begin{array}{cc}
        -\left( -\frac{I}{2} +  \td{\bK}_{\partial D}^{\omega, *} \right)^{-1} \left( \frac{I}{2} +  {\bK}_{\partial D}^{\omega, *} \right) \left( {\bS}_{\partial D}^{\omega} \right)^{-1} &  \left( -\frac{I}{2} +  \td{\bK}_{\partial D}^{\omega, *} \right)^{-1}\medskip \\
      -\left( {\bS}_{\partial D}^{\omega} \right)^{-1} & 0 \\
    \end{array}
  \right).
\end{split}
\]
Moreover, through tedious calculation and from the asymptotic expansion for $\omega\ll1$ and $\delta\ll1$, one has that 
\[
\norm{\Gcal_1^{-1}(\omega,\delta)\Gcal_2(\omega,\delta)}_{\Lcal(\Hcal,\Hcal)} = \Ocal\left( \frac{\epsilon}{\omega^2}\right).
\]
From Lemma \ref{lem:sge}, one can conclude that there exist characteristic values $\omega^*\ll1$, located in the region $V$, of the operator $\Acal(\omega,\delta)$. This completes the proof.

% Consider the kernel of the operator $\Gcal(\omega,\delta)$, i.e.
% \begin{equation}\label{eq:keg}
%     \Gcal(\omega,\delta)  \left(
%     \begin{array}{c}
%        \mathbf{g}_1 \\
%       \mathbf{g}_2\\
%     \end{array}
%   \right).
% \end{equation}
% The first equation in \eqref{eq:keg} is equivalent to 
% \[
% {\bS}_{\partial D}^{\omega}[\mathbf{g}_2](\bx)=0, \quad \bx\in\partial D.
% \]
% If $\mathbf{g}_2$ is nontrivial, the function defined by
% \[
% \bu(\bx)={\bS}_{\partial D}^{\omega}[\mathbf{g}_2](\bx), \quad \bx\in \mathbb{R}^3\backslash {D},
% \]
% will solve the following Dirichlet boundary problem
% \begin{equation}
%   \left\{
%     \begin{array}{ll}
%        \mathcal{L}_{\lambda, \mu}\bu(\bx) + \omega^2\rho\bu(\bx) =0    &   \bx\in \mathbb{R}^3\backslash \overline{D},  \medskip \\
%       \bu(\bx) = 0      & \bx\in\partial D,  \medskip \\
%       \bu  \quad \mbox{satisfies the radiation condition}.
%     \end{array}
%   \right.
% \end{equation} 
% The above problem only has the trivial solution [...], which shows that the operator ${\bS}_{\partial D}^{\omega}$ is invertible from $L^2(\partial D)^3$ to $H^1(\partial D)^3$.

\end{proof}

After obtaining the existence result for the sub-wavelength resonances, we then explicitly derive the mathematical formula for the resonant frequencies. 
For the further analysis, we define the projection operator $\Pcal_i$ as
\begin{equation}\label{eq:dp}
 \Pcal_i[\Phi] = (\Phi, \hat\Phi_i) \hat\Psi_i, \quad 1\leq i \leq 6,
\end{equation}
where functions $\hat{\Phi}_i$ and $\hat\Psi_i$ are given in Lemma \ref{lem:ker}.
We also define the operator 
\begin{equation}\label{eq:hao}
 \hat\Acal_0:= \Acal_0 + \Pcal = \Acal_0 +\sum_{i=1}^{6} \Pcal_i,
\end{equation}
where the operator $\Acal_0$ is given in Lemma \ref{lem:deoa}.
Then the operator $ \hat\Acal_0$ shares the following properties.

\begin{lem}
We have
\begin{enumerate}
\item the operator $\hat\Acal_0$ defined in \eqref{eq:hao} is bijective in $\Lcal(\Hcal,\Hcal^1)$ and $\hat\Acal_0[\hat\Phi_i] =  \hat\Psi_i$, with $1\leq i \leq 6$;
\item the adjoint operator of  $\hat\Acal_0$, i.e. $\hat\Acal_0^*$, is also bijective and $\hat\Acal_0^*[\hat\Psi_i] =  \hat\Phi_i$.
\end{enumerate}
\end{lem}

\begin{proof}
Note that the operator $\Acal_0$ is a Fredholm operator with the index $0$. Thus the bijection of the operator $\hat\Acal_0$ directly follows from that the construction of the operator $\hat\Acal_0$ in \eqref{eq:hao}. Direct calculation shows that 
\[
 \hat\Acal_0[\hat\Phi_i]= \Acal_0[\hat\Phi_i] +\sum_{j=1}^{6} \Pcal_j [\hat\Phi_i]=\delta_{ij} \hat\Psi_j = \hat\Psi_i.
\]
Thus the first statement is proved. Next, we prove the second one. Direct calculation shows that 
\[
 \hat\Acal_0^*[\hat\Psi_i]  =  \Acal_0^*[\hat\Psi_i]  + \Pcal^*[\hat\Psi_i]  = \sum_{j=1}^{6} \Pcal_j^* [\hat\Phi_i]=\hat\Phi_i
\]
This completes the proof.
\end{proof}

After these preparations, we are in a position to derive the formula of the sub-wavelength resonant frequencies for the system \eqref{eq:mo}, which is concluded in the following theorem. 
 
\begin{thm}\label{thm:ref}
Consider the system \eqref{eq:mo} with the parameters given in \eqref{eq:hicon} and \eqref{eq:de}. Then the sub-wavelength resonant frequencies have the following expressions, 
\[
\omega_i = \epsilon^{1/2} \left( \frac{1}{\rho\varrho_i } \right)^{1/2} (1+ o(1)),
\]
where $\rho$ is the density in the background, $\epsilon$ is defined in \eqref{eq:hicon} and $\varrho_i, 1\leq i\leq 6$ are the eigenvalues of the matrix $M$ defined in Lemma \ref{lem:ma1}.
\end{thm}
\begin{proof}
We would like to find the characteristic values $\omega^*\ll1$ such that there exits a nontrivial function $\Phi_\delta$ satisfying 
\begin{equation}\label{eq:eqa}
\Acal(\omega^*, \delta)[\Phi_\delta]=0.
\end{equation}
Lemma \ref{lem:ker} shows that 
\[
\Acal(0,0)[\Phi_0] =0,
\]
where $\Phi_0 =\sum_{j=1}^6 c_j \hat\Phi_j$ and the coefficients $c_j$ are arbitrary. 
Thus one can treat $\Phi_\delta$ as a perturbation of $\Phi_0$. We express the function $\Phi_\delta$ by 
\[
 \Phi_\delta = \Phi_0 + \Phi_1
\]
where the function $\Phi_1$ is uniquely determined by requiring 
\begin{equation}\label{eq:orso}
 (\Phi_1, \Phi_0)_{\Hcal}=0.
\end{equation}
Using the operator $\hat\Acal_0$ defined in \eqref{eq:hao}, the equation \eqref{eq:eqa} is equivalent to 
\begin{equation}\label{eq:exa0}
 \left( \hat\Acal_0 + \Bcal -\Pcal  \right)[\Phi_0 + \Phi_1]  =0.
\end{equation}
From Lemma \ref{lem:deoa}, the operator $ \hat\Acal_0 + \Bcal$ is invertible for sufficient small $\omega$ and $\delta$. Moreover, due to $\norm{\Bcal}_{\Lcal(\Hcal,\Hcal)}\ll 1$, the operator $ (\hat\Acal_0 + \Bcal)^{-1}$ has the following asymptotic expansion
\begin{equation}
\begin{split}
  &(\hat\Acal_0 + \Bcal)^{-1} =  (I +  \hat\Acal_0^{-1} \Bcal)^{-1}  \hat\Acal_0^{-1} \\
 = &  \left( I  - \hat\Acal_0^{-1} \Bcal + ( \hat\Acal_0^{-1} \Bcal)^2 + \cdots \right)   \hat\Acal_0^{-1}.
\end{split}
\end{equation}
Then applying $(\hat\Acal_0 + \Bcal)^{-1}$ on both sides of the equation \eqref{eq:exa0} yields that 
\[
\Phi_1 = (\hat\Acal_0 + \Bcal)^{-1}\Pcal[\Phi_0] - \Phi_0.
\]
Thus the equation \eqref{eq:orso} can be expanded as  
\begin{equation}\label{eq:incal}
\begin{split}
 (\Phi_1, \Phi_0)_{\Hcal} =& \left(  (\hat\Acal_0 + \Bcal)^{-1}\Pcal[\Phi_0] - \Phi_0, \Phi_0  \right)_{\Hcal}   \\
%= & \left(  (\hat\Acal_0 + \Bcal)^{-1} \Psi_0 , \Phi_0  \right)_{\Hcal}   -   \left( \Phi_0, \Phi_0  \right)_{\Hcal}  \\
= & \left(  \left( I  - \hat\Acal_0^{-1} \Bcal + ( \hat\Acal_0^{-1} \Bcal)^2 + \cdots \right)   \Phi_0 , \Phi_0  \right)_{\Hcal}   -   \left( \Phi_0, \Phi_0  \right)_{\Hcal}  \\
= & \left(  \left(  - \hat\Acal_0^{-1} \Bcal + ( \hat\Acal_0^{-1} \Bcal)^2 + \cdots \right)   \Phi_0 , \Phi_0  \right)_{\Hcal} .   \\
\end{split}
\end{equation}
In the derivation of \eqref{eq:incal}, we have used the following fact
\[
\Pcal[\Phi_0] = \Psi_0 = \sum_{j=1}^6 c_j  \hat\Psi_j,
\]
which follows from the definition of the operator $\Pcal$ in \eqref{eq:dp}.
Next, we do more analysis on the equation \eqref{eq:incal}.
From the expressions of operators $\Acal_{m,0}$ in \eqref{lem:deoa}, the following identities hold for $m\geq 1$
\[
\Acal_{m,0}  \Phi_0 =0. 
\]
Thus we have that 
\begin{equation}\label{eq:esi2}
 ( \hat\Acal_0^{-1} \Bcal)^2  \Phi_0 =  \Ocal \left( \delta\omega + \delta^2 + \omega^2\tau^2\delta  \right). 
\end{equation}
Then we consider this term $ \left(   - \hat\Acal_0^{-1} \Bcal   \Phi_0 , \Phi_0  \right)_{\Hcal} $.
Direct calculation shows that 
\begin{equation}\label{eq:incal2}
\begin{split}
  &\left(   - \hat\Acal_0^{-1} \Bcal   \Phi_0 , \Phi_0  \right)_{\Hcal} 
=  \left(   - \Bcal   \Phi_0 , ( \hat\Acal_0^{-1})^* \Phi_0  \right)_{\Hcal^1} 
=  \left(   - \Bcal   \Phi_0 , \Psi_0  \right)_{\Hcal^1}   \\
= &  \left(   - \left( \delta\Acal_{0,1} + \omega^2 \tau^2  \Acal_{0,2}  \right)   \Phi_0 , \Psi_0  \right)_{\Hcal^1} + \Ocal \left(  \omega \delta \tau + \left( \omega \tau \right)^3  \right ) \\
 = & \sum_{i,j=1}^6 c_i c_j d_j\left( \delta \left( \bS_{\partial D} [\bzeta_i] , (\bS^{-1}_{\partial D})^*[ \bxi_j] \right)_{H^1} - \omega^2 \tau^2 \left(\bK_{\partial D, 2}^{*}[\bzeta_i] ,   \bxi_j  \right)_{L^2} \right) +  \Ocal \left(  \omega \delta \tau + \left( \omega \tau \right)^3  \right )  \\
 = & \sum_{i,j=1}^6 c_i c_j d_j\left( \delta \left( \bzeta_i ,  \bxi_j \right)_{L^2}  - \omega^2 \tau^2  \left(\bK_{\partial D, 2}^{*}[\bzeta_i] ,   \bxi_j  \right)_{L^2} \right) +   \Ocal \left(  \omega \delta \tau + \left( \omega \tau \right)^3  \right ). 
\end{split}
\end{equation}

Furthermore, using Green's formula yields that
\begin{equation}\label{eq:ink2}
\begin{split}
& \left(\bK_{\partial D, 2}^{*}[\bzeta_i] ,   \bxi_j  \right)_{L^2}  =  \int_{\partial D} \bK_{\partial D, 2}^{*}[\bzeta_i](\bx) \cdot   \bxi_j(\bx) ds(\bx)  =  \int_{\partial D} \bzeta_i(\bx) \cdot    \bK_{\partial D, 2}[\bxi_j](\bx) ds(\bx) \\
= &- \rho \int_{\partial D} \bzeta_i(\bx) \cdot   \int_D    \mathbf{\Gamma}_0(\bx-\by) \bxi_j(\by) d\by ds(\bx) =- \rho \int_{D} \bxi_j(\by) \cdot   \int_{\partial D}    \mathbf{\Gamma}_0(\bx - \by)\bzeta_i(\bx) ds(\bx) d\by \\
= & - \rho \int_{D} \bxi_j(\by) \cdot  \bS_{\partial D}[\bzeta_i](\by) d\by.
\end{split}
\end{equation}
%In the derivation of the last equation, we have used the following identity 
%\[
%  \mathcal{L}_{\lambda, \mu}  \mathbf{\Gamma}_2(\bx) + \rho  \mathbf{\Gamma}_0(\bx)  =0,
%\]
%where $\mathbf{\Gamma}_0$ and $\mathbf{\Gamma}_2$ are given in \eqref{eq:fse}.
Thus using Lemma \ref{lem:orker} and combining the equations \eqref{eq:incal} -- \eqref{eq:ink2}, we finally have that 
\[
\begin{split}
(\Phi_1, \Phi_0)_{\Hcal} = & \sum_{i,j=1}^6 c_i c_j d_j\left( \delta \delta_{ij} - \omega^2 \tau^2 M_{ij}\right)+  \Ocal \left(  \omega \delta \tau + \left( \omega \tau \right)^3  \right ) \\
= &  \bc^t \left( \delta I  - \omega^2 \tau^2 M  \right) \Lambda\bc + \Ocal \left(  \omega \delta \tau + \left( \omega \tau \right)^3  \right )=0,
\end{split}
\]
where $\bc=(c_1, c_2, c_3, c_4, c_5, c_6)^t$, $\Lambda$ is a diagonal matrix with $\Lambda_{ii}=d_i$ and $M$ is a matrix defined in \eqref{eq:mij}.
Thus the resonant frequencies are 
\[
 \omega_i = \epsilon^{1/2} \left( \frac{1}{\rho\varrho_i } \right)^{1/2} (1+ o(1)),
\]
where parameters $\varrho_i, 1\leq i\leq 6$ are the eigenvalues of the matrix $M$.
This completes the proof.
\end{proof}

\begin{rem}\label{rem:depre}
Here we would like to remark that the resonant phenomenon is caused by the high contrast of the Lam\'e parameters in different region, which can be concluded from Lemma \ref{lem:exres}. Indeed, from the expression of the operator $\Acal(\omega,\delta)$ in Lemma \ref{lem:deoa} and the fact that operator $\Acal_0$ has a nontrivial kernel from Lemma \ref{lem:ker}, the norm of the operator $ \td{\bS}_{\partial D}^{\omega}$ should be small; that is $\norm{\td{\bS}_{\partial D}^{\omega}}_{\Lcal(\Hcal,\Hcal^1)}=o(1)$. Moreover, direct calculation shows that $\norm{\td{\bS}_{\partial D}^{\omega}}_{\Lcal(\Hcal,\Hcal^1)}=\Ocal(\delta)$. Thus to ensure the occurrence of the resonances, 
 the contrast of the Lam\'e parameters in different region should be high. Thus one can conclude that the resonance is caused by the high  contrast of the Lam\'e parameters. Furthermore, from the expression of the resonant frequencies in Theorem \ref{thm:ref}, one can have that the sub-wavelength characteristic of the resonances is determined by the high contrast of the mass densities in different regions. 
\end{rem}

\subsection{Dipolar resonances}

In this subsection, we study the behavior of the displacement filed inside the domain $D$ for the system \eqref{eq:mo} when the incident frequency $\omega$ is located in different regimes. In what follows, we take the incident wave $\bu^i$ in \eqref{eq:mo} to be a p-wave, namely
\begin{equation}\label{eq:inp}
\bu^i= \bp e^{-\rmi k_p \bp\cdot\bx}  ,
\end{equation}
where $\bp\in\mathbb{R}^3$ denotes the polariton direction and satisfies $\bp\cdot\bp=1$. For other incident waves, one can obtain similar conclusions following the same discussion.
We have the following theorem about the displacement filed $\bu_D$ inside the domain $D$ for the frequency $\omega$ located in different regimes.
\begin{thm}\label{thm:dir}
Consider the system \eqref{eq:mo} with the parameters given in \eqref{eq:hicon} and \eqref{eq:de}. Then the displacement filed $\bu_D$ inside the domain $D$ has the following expression for $\omega\ll1$:
\begin{equation}
\bu_D = \sum_{i=1}^{6} \gamma_i \bxi_i (1 + o(1)),
\end{equation}
where $\bxi_i$, $1\leq i\leq 6$, are given in Lemma \ref{lem:orker}.
Moreover, when the frequency $\omega$ is located in the following different regimes, the coefficients $\gamma_i, 1\leq i\leq 6$, satisfy
\begin{enumerate}
\item for $\omega\ll \sqrt{\epsilon}$, the coefficients $\gamma_i=\Ocal(1)$;
\item for $\omega = \Ocal( \sqrt{\epsilon})$, the coefficients $\gamma_i=\Ocal \left( \delta^{-1/2} \right )$;
\item for $ \sqrt{\epsilon} \ll \omega \ll1$, the coefficients $\gamma_i= o(1)$.
\end{enumerate}
\end{thm}

\begin{proof}
From the solution written in \eqref{eq:sol}, to investigating the system \eqref{eq:mo} we need to solve the equation \eqref{eq:or}. With the help of the operator $\Pcal$ defined in \eqref{eq:hao}, the equation \eqref{eq:or} is equivalent to 
\begin{equation}\label{eq:orf}
\left( \hat\Acal_0 + \Bcal -\Pcal  \right)[\Phi](\bx) = F(\bx)  \quad \bx\in\partial D,
\end{equation}
where $\Phi\in\Hcal$. Then we focus ourself on the equation \eqref{eq:orf}.
We first decompose the function $\Phi$ into two parts, i.e. 
\[
 \Phi=\Phi_k + \Phi_t,
\]
where $\Phi_k\in \ker \Acal_{0}$ with $\ker \Acal_{0}$ defined in Lemma \ref{lem:ker} and $\Phi_t$ is uniquely determined by $(\Phi_k, \Phi_t)_{\Hcal}=0$. Thus the function $\Phi_k$ can be written as 
\begin{equation}\label{eq:pk1}
 \Phi_k = \sum_{i=1}^6 \alpha_i \hat\Phi_i,
\end{equation}
where the coefficients $\alpha_i$ shall be determined later. Note that the operator $\hat\Acal_0 + \Bcal$ is invertible due to $\norm{\Bcal}_{\Lcal(\Hcal,\Hcal^1)}\ll 1$ and the invertibility of the operator $\hat\Acal_0 $.
Multiplying $(\hat\Acal_0 + \Bcal)^{-1}$ on both sides of equation \eqref{eq:orf} yields that
\begin{equation}
\left( I - (\hat\Acal_0 + \Bcal)^{-1}\Pcal \right)[\Phi_k + \Phi_t] = (\hat\Acal_0 + \Bcal)^{-1}[F],
\end{equation}
and further simplification gives that 
\begin{equation}\label{eq:simor}
\Phi_k + \Phi_t - \left( I +  \hat\Acal_0^{-1} \Bcal \right)^{-1} [\Phi_k] = (\hat\Acal_0 + \Bcal)^{-1}[F].
\end{equation}
Next, we analyze the right hand side of equation \eqref{eq:simor}, i.e. $(\hat\Acal_0 + \Bcal)^{-1}[F]$. If the incident wave $\bu^i$ is given in \eqref{eq:inp}, the function $F$ in \eqref{eq:simor} has the following expansion
\begin{equation}\label{eq:decF}
F= F_0 + \Ocal(\omega),
\end{equation} 
where 
\begin{equation}
 F_0= \left(
    \begin{array}{c}
     \bp \\
     0 \\
    \end{array}
  \right).
\end{equation}
%Thus the polariton direction $\bp$ can be expressed by 
%\[
% \bp = \sum_{i=1}^{6}  f_i \bxi_i,
%\]
%for certain coefficients $f_i=(\bp, \bxi_i)$. 
%Denote by 
%\begin{equation}
% F_0= \left(
%    \begin{array}{c}
%     \bp \\
%     0 \\
%    \end{array}
%  \right),
%\end{equation}
%and 
Further calculation shows that $F_0$ enjoys the following decomposition,
\begin{equation}
F_0 = F_1 + \sum_{i=1}^{6} (F_0,  \hat\Psi_i)_{\Hcal^1}  \hat\Psi_i,
\end{equation}
with $\beta_i=(F_0,  \hat\Psi_i)_{\Hcal^1}$, the functions $\hat\Psi_i$ given in Lemma \ref{lem:ker} and $F_1 = F_0 - \sum_{i=1}^{6} (F_0,  \hat\Psi_i)_{\Hcal^1}  \hat\Psi_i$.
By the definition of the operator $\hat\Acal_0$ in \eqref{eq:hao}, the following holds 
\begin{equation}\label{eq:inv2}
\hat\Acal_0^{-1} \left[F_0 \right] = \Acal_0^{-1}[F_1] + \sum_{i=1}^{6} \beta_i \hat\Phi_i,
\end{equation}
where 
\begin{equation}
\left( \Acal_0^{-1}[F_1] \right)_1 \notin \ker \left( -1/2 +  \bK^*_{\partial D} \right),
\end{equation}
with $\left( \Acal_0^{-1}[F_1] \right)_1 $ signifying the first component of the vector $ \Acal_0^{-1}[F_1] $. Finally, combining equations from \eqref{eq:decF} to \eqref{eq:inv2}, the right hand side of \eqref{eq:simor} can be written as 
\begin{equation}\label{eq:rf}
\begin{split}
 &(\hat\Acal_0 + \Bcal)^{-1}[F] = (\hat\Acal_0 + \Bcal)^{-1}\left[F_0 \right] + \Ocal(\omega )\\
  = & \Acal_0^{-1}[F_1] + \sum_{i=1}^{6} \beta_i \hat\Phi_i + \Ocal(\omega + \delta).
\end{split}
\end{equation}
%
%Then one has that 
%\begin{equation}\label{eq:inv1}
%(\hat\Acal_0 + \Bcal)^{-1}[F] = (\hat\Acal_0 + \Bcal)^{-1}\left[F_0 \right] + \Ocal(\omega ).
%\end{equation}
%
%Combining equations \eqref{eq:inv1} and \eqref{eq:inv2}, one finally has that 
%\begin{equation}
%(\hat\Acal_0 + \Bcal)^{-1}[F] =  \Acal_0^{-1}[F_1] + \sum_{i=1}^{6} (F_0,  \hat\Psi_i)_{\Hcal^1}  \hat\Phi_i + \Ocal(\omega + \delta).
%\end{equation}
Then, we determinate the coefficients $\alpha_i$ in \eqref{eq:pk1}. Multiplying $\hat\Phi_i$ on both sides of \eqref{eq:simor}, integrating on $\partial D$ and together with the help of \eqref{eq:rf}, one has that the coefficients $\alpha_i$ in \eqref{eq:pk1} satisfy
\begin{equation}\label{eq:eqma}
 \Lambda \left( \delta I  - \omega^2 \tau^2 M^t  \right) \bal + \Ocal \left(  \omega \delta \tau + \left( \omega \tau \right)^3  \right ) = \bbe +  \Ocal(\omega + \delta),
\end{equation}
where $\bal=(\alpha_1, \alpha_2, \alpha_3, \alpha_4, \alpha_5, \alpha_6)^t$,  $\bbe=(\beta_1, \beta_2, \beta_3, \beta_4, \beta_5, \beta_6)^t$, $\Lambda$ is a diagonal matrix with $\Lambda_{ii}=d_i$ and $M$ is a matrix defined in \eqref{eq:mij}. 
We would like to mention that the order of $ \Ocal \left(  \omega \delta \tau + \left( \omega \tau \right)^3  \right ) $ on the left hand side of \eqref{eq:eqma} is imaginary from the operator expansion in Lemma \ref{lem:deoa}.
%Thus, the coefficients $\alpha_i$ in \eqref{eq:pk1} satisfy
%\begin{equation}
% \bal =  \left( \Lambda \left( \delta I  - \omega^2 \tau^2 M^t  \right) \right)^{-1}\bbe + o(1). 
%\end{equation}
From \eqref{eq:simor} and \eqref{eq:decF}, one has that 
\[
  \Phi_t = \hat\Acal_0^{-1} \left[F_0 \right] +   \sum_{i=1}^{6} \alpha_i \left( \delta \hat\Acal_0^{-1} \Acal_{0,1} [ \hat\Phi_i ]+ \omega^2 \tau^2 \Acal_0^{-1} \Acal_{0,2} [ \hat\Phi_i ]  \right)  +o(1).
\]
Finally, the displacement field inside the domain $D$ is given by 
\begin{equation}
\begin{split}
 &\bu_D= \td{\bS}_{\partial D}^{\omega}[\Phi_k + \Phi_t] =  \delta\bS_{\partial D}[\Phi_k + \Phi_t] + \Ocal(\delta\omega) =  \delta \left( \bS_{\partial D}[\Phi_k ] + \bS_{\partial D}[ \Phi_t] \right) + \Ocal(\delta\omega) \\
  = & \delta \left(  \sum_{i=1}^6 \alpha_i \bS_{\partial D} [  \bzeta_i ] \left( 1 + \Ocal \left(\delta + \omega^2\tau^2 \right) \right)  +  \bS_{\partial D} \left[  \hat\Acal_0^{-1} \left[F_0 \right]  \right]  \right) + \Ocal(\delta\omega).
\end{split}
\end{equation}
Thus, when $\omega\ll \sqrt{\delta}/\tau= \sqrt{\epsilon}$, one has that $\alpha=\Ocal(\delta^{-1})$ from \eqref{eq:eqma} and displacement field inside the domain $D$ statisfies
\[
 \bu_D =\delta \sum_{i=1}^6 \alpha_i \bS_{\partial D} [  \bzeta_i ] + \Ocal(\delta) = \Ocal(1).
 \]
%where the leading term is of order $1$; that is 
%\[
%\delta \sum_{i=1}^6 \alpha_i \bS_{\partial D} [  \bzeta_i ] = \Ocal(1).
%\]
When $\omega = \Ocal( \sqrt{\delta}/\tau)=\Ocal( \sqrt{\epsilon})$,  one has that 
%\[
% \alpha=   \frac{1} {\delta \left(  I  - \omega^2 \tau^2/\delta M^t + \Ocal \left(  \omega  \tau  \right )   \right)   }   ,
%\]
\[
 \alpha=   \frac{1} {\delta   \Ocal \left( \sqrt{\delta} \right )    }   ,
\]
and displacement field inside the domain $D$ satisfies
\[
 \bu_D =\delta \sum_{i=1}^6 \alpha_i \bS_{\partial D} [  \bzeta_i ] + \Ocal(\delta^{1/2}) =  \Ocal \left( \delta^{-1/2} \right ).
 \]
%where the leading term is of order $\Ocal \left( \delta^{-1/2} \right ) $; that is 
%\[
%\delta \sum_{i=1}^6 \alpha_i \bS_{\partial D} [  \bzeta_i ] = \Ocal \left( \delta^{-1/2} \right ) .
%\]
When $   \sqrt{\epsilon} =  \sqrt{\delta}/\tau \ll \omega \ll1$,  one has that 
%\[
% \alpha=   \frac{1} {\delta \left(  I  - \omega^2 \tau^2/\delta M^t + \Ocal \left(  \omega  \tau  \right )   \right)   }   ,
%\]
\[
 \alpha=   \frac{1} {   \Ocal \left( \omega^2\tau^2 \right )    }   ,
\]
and displacement field inside the domain $D$
\[
 \bu_D =\delta \sum_{i=1}^6 \alpha_i \bS_{\partial D} [  \bzeta_i ] + \Ocal(\delta)= \Ocal \left( \delta/(\omega\tau)^2 \right )=o(1).
 \]
%where the leading term is of order $\Ocal \left( \delta/(\omega\tau)^2 \right ) $; that is 
%\[
%\delta \sum_{i=1}^6 \alpha_i \bS_{\partial D} [  \bzeta_i ] = \Ocal \left( \delta/(\omega\tau)^2 \right )=o(1) .
%\]
The proof is completed by noting that the operator $\bS_{\partial D} $ is bounded from $\mho^*$ to $\mho$.
\end{proof}

\begin{rem}
From Remark \ref{re:dipb}, the fields $\bxi_i$, $1\leq i\leq 6$ can be treated as dipolar wave fields. The last theorem shows that if the frequency $\omega$ is of the order $\Ocal(\sqrt{\epsilon} )$, i.e. around the resonant frequencies from Theorem \ref{thm:ref}, the dipolar wave fields inside the domain $D$ is greatly enhanced at the order $\Ocal \left( \delta^{-1/2} \right )$.
\end{rem}

\section {The analysis of the resonances for the spherical geometry}

In this section, we discuss the sub-wavelength resonant phenomenon for the system \eqref{eq:mo} when the domain $D$ is a ball. We apply two methods to study the resonant phenomenon. One is to use the conclusions obtained in Theorem \ref{thm:ref}. The other is to directly solve the equation \eqref{eq:conre1} to obtain the resonant frequencies. At last, the characterization of the dipolar resonance is presented. For the convenience of readers, we first present some preliminaries about the properties of the single-layer potential operators and the N-P operators in the spherical geometry. 

Let $j_n(t)$ and $h_n(t)$, $n\in\mathbb{N}_0$, denote the spherical Bessel and Hankel functions of the first kind with the order $n$, respectively.
The following asymptotic expansions hold for $t\ll 1$ \cite{CK3518},
\begin{equation}\label{eq:asjt}
\begin{split}
  j_1(t)=  \frac{t}{3}  -  \frac{t^3}{30}  + \Ocal(t^5), \quad
  h_1(t)=-\frac{\ri}{t^2}  - \frac{\ri}{2} + \frac{t}{3} +  \frac{\ri t^2}{8}  + \Ocal(t^3),
 \end{split}
\end{equation}
Let $Y_n^m(\hat{\bx})$ denote spherical harmonic functions of the order $n$ with the degree $m$. Denote by $B_R$ a centred ball with the radius $R$. 
Let $\bbS_{R}$ be the surface of the ball $B_R$ and denote by $\bbS$ for $R=1$ for simplicity. 
The operator $\nabla_{\bbS}$ designates the surface gradient on the unit sphere $\bbS$.
Then the following lemma holds \cite{DLL7678, L0069}.
\begin{lem}\label{lem:vec}
The family $(\Ical_n^m, \Tcal_n^m, \Ncal_n^m)$, the vectorial spherical harmonics of order $n$,
\[
\begin{split}
 \Ical_n^m= &\nabla_{\bbS}Y_{n+1}^m +(n+1) Y_{n+1}^m \bnu, \quad n\geq0, \; n+1\geq m \geq -(n+1),\\
 \Tcal_n^m= &\nabla_{\bbS}Y_{n}^m\wedge\bnu, \qquad\qquad\qquad  n\geq1, \; n\geq m \geq -n,\\
 \Ncal_n^m=&-\nabla_{\bbS}Y_{n-1}^m + nY_{n-1}^m \bnu, \quad n\geq1, \; n-1\geq m \geq -(n-1),
\end{split}
\]
forms an orthogonal basis of $(L^2(\bbS))^3$, where $\bnu$ denotes the outward unit normal to unit sphere.
\end{lem}
From the definition of the vectorial spherical harmonics in Lemma \ref{lem:vec}, direct calculation yields that 
\begin{equation}\label{eq:I}
  \Ical_0^{-1}=\sqrt{\frac{3}{8\pi}}(1, -\ri, 0)^t, \quad  \Ical_0^{0}=\sqrt{\frac{3}{4\pi}}(0, 0, 1)^t, \quad   \Ical_0^{1}=-\sqrt{\frac{3}{8\pi}}(1, \ri, 0)^t, 
\end{equation}
and 
\begin{equation}\label{eq:T}
\begin{split}
 r\Tcal_1^{-1}=\sqrt{\frac{3}{8\pi}} & (-\ri x_3, -x_3, \ri x_1 + x_2)^t, \quad  r\Tcal_1^{0}=\sqrt{\frac{3}{4\pi}}(-x_2, x_1, 1)^t, \\
  & r\Tcal_1^{1}=-\sqrt{\frac{3}{8\pi}}(\ri x_3, -x_3, -\ri x_1 + x_2)^t.
\end{split}
\end{equation}
Thus the vectorial spherical harmonics $\Ical_0^m$ and $r\Tcal_1^m$, $-1\leq m\leq 1$, forms a basis for the space $\mho$ given in \eqref{eq:dpsi}. For the single-layer potential and the N-P operators, one has the following two lemmas from \cite{DLL7678, L0069}.

\begin{prop}\label{pro:singlei}
For $\bx\in B_R$, the single-layer potentials $ \bS_{\bbS_R}^{\omega}[\Tcal_n^m](\bx)$, $ \bS_{\bbS_R}^{\omega}[\Ical_{n-1}^m](\bx)$ and $  \bS_{\bbS_R}^{\omega}[\Ncal_{n+1}^m](\bx)$  are given as follows 
\[
\begin{split}
  \bS_{\bbS_R}^{\omega}[ \Ical_{n-1}^m](\bx) =&  -R^2 \rmi \left( \frac{ (n+1)k_s h_{n-1,s} j_{n-1}(k_s |\bx|) }{\mu(2n+1)}  +  \frac{ nk_p h_{n-1,p} j_{n-1}(k_p |\bx|)}{(\lambda+2\mu)(2n+1)}   \right) \Ical_{n-1}^m \\
       &- n R^2\rmi  \left( \frac{k_s h_{n-1,s} j_{n+1}(k_s |\bx|) }{\mu(2n+1)}  -   \frac{k_p h_{n-1,p} j_{n+1}(k_p |\bx|) }{(\lambda+2\mu)(2n+1)}  \right) \Ncal_{n+1}^m,
 \end{split}
\]

\[
\begin{split}
  \bS_{\bbS_R}^{\omega}[ \Ncal_{n+1}^m](\bx) =&  - (n+1) R^2\rmi  \left( \frac{ k_s h_{n+1,s} j_{n-1}(k_s |\bx|)}{\mu(2n+1)}  -  \frac{ k_p h_{n+1,p} j_{n-1}(k_p |\bx|) }{(\lambda+2\mu)(2n+1)}   \right) \Ical_{n-1}^m \\
       &-R^2 \rmi \left( \frac{ n k_s h_{n+1,s} j_{n+1}(k_s |\bx|) }{\mu(2n+1)}  +  \frac{ (n+1) k_p h_{n+1,p} j_{n+1}(k_p |\bx|)}{(\lambda+2\mu)(2n+1)}  \right) \Ncal_{n+1}^m,
 \end{split}
\]
and 
\[
 \bS_{\bbS_R}^{\omega}[\Tcal_n^m](\bx) = -\frac{\rmi k_s R^2 h_{n,s} j_n(k_s|\bx|)}{\mu}\Tcal_n^m.
\]

\end{prop}

\begin{lem}\label{lem:sin}
For $\bx\in \bbS_R$, the single-layer potentials $ \bS_{\bbS_R}^{\omega}[\Tcal_n^m]$, $ \bS_{\bbS_R}^{\omega}[\Ical_{n-1}^m]$ and $  \bS_{\bbS_R}^{\omega}[\Ncal_{n+1}^m]$  are given as follows 
\begin{equation*}\label{eq:sint}
\bS_{\bbS_R}^{\omega}[\Tcal_n^m](\bx) = b_n\Tcal_n^m,\quad
 \bS_{\bbS_R}^{\omega}[ \Ical_{n-1}^m](\bx) = c_{1n} \Ical_{n-1}^m + d_{1n} \Ncal_{n+1}^m,
\end{equation*}
and
\begin{equation*}
 \bS_{\bbS_R}^{\omega}[ \Ncal_{n+1}^m](\bx) = c_{2n} \Ical_{n-1}^m + d_{2n}\Ncal_{n+1}^m,
\end{equation*}
where
\[
 b_n=  -\frac{\rmi k_s R^2 j_{n,s} h_{n,s}}{\mu},
\]
\[
\begin{split}
c_{1n}=& -R^2 \rmi \left( \frac{j_{n-1,s} h_{n-1,s} k_s (n+1)}{\mu(2n+1)}  +  \frac{j_{n-1,p} h_{n-1,p} k_p n}{(\lambda+2\mu)(2n+1)}   \right),\\
d_{1n}=& - n R^2\rmi  \left( \frac{j_{n-1,s} h_{n+1,s} k_s}{\mu(2n+1)}  -   \frac{j_{n-1,p} h_{n+1,p} k_p }{(\lambda+2\mu)(2n+1)}  \right),
\end{split}
\]
\[
\begin{split}
c_{2n}=& - (n+1) R^2\rmi  \left( \frac{j_{n+1,s} h_{n-1,s} k_s}{\mu(2n+1)}  -  \frac{j_{n+1,p} h_{n-1,p} k_p }{(\lambda+2\mu)(2n+1)}   \right),\\
d_{2n}=&  -R^2 \rmi \left( \frac{j_{n+1,s} h_{n+1,s} k_s n}{\mu(2n+1)}  +  \frac{j_{n+1,p} h_{n+1,p} k_p (n+1)}{(\lambda+2\mu)(2n+1)}  \right).
\end{split}
\]

\end{lem}

\begin{lem}\label{lem:trac}
For the N-P operator $\bK_{\partial D}^{\omega, *}$ with the density functions $(\Ical_n^m, \Tcal_n^m, \Ncal_n^m)$, the following identities hold on $\bbS_R$
\begin{eqnarray*}
\bK_{\partial D}^{\omega, *}[\Tcal_n^m](\bx) & = & (\mathfrak{b}_n -1/2) \Tcal_n^m,\label{eq:tract}\\
\bK_{\partial D}^{\omega, *}[ \Ical_{n-1}^m] (\bx) & = & (\mathfrak{c}_{1n} - 1/2) \Ical_{n-1}^m +  \mathfrak{d}_{1n} \Ncal_{n+1}^m,\label{eq:traci}\\
\bK_{\partial D}^{\omega, *}[ \Ncal_{n+1}^m] (\bx) & = & \mathfrak{c}_{2n} \Ical_{n-1}^m +  (\mathfrak{d}_{2n}-1/2) \Ncal_{n+1}^m,\label{eq:tracn}
\end{eqnarray*}
where
\[
  \mathfrak{b}_n=  -\rmi k_s R j_{n,s} (k_sR h_{n,s}^{\prime} - h_{n,s}),
\]
\[
 \begin{split}
 \mathfrak{c}_{1n}=& -2(n-1)R\rmi  \left( \frac{j_{n-1}(k_s R) h_{n-1}(k_s R) k_s (n+1)}{2n+1}  +  \frac{j_{n-1}(k_p R) h_{n-1}(k_p R) k_p\mu n}{(\lambda+2\mu)(2n+1)}   \right)\\
         & + R^2\rmi  \left( \frac{ j_{n-1}(k_s R) h_{n}(k_s R) k_s^2  (n+1) + j_{n-1}(k_p R) h_{n}(k_p R) k_p^2 n}{2n+1}   \right), \\
 \mathfrak{d}_{1n}=&  2n(n+2) R\rmi  \left( \frac{j_{n-1}(k_s R) h_{n+1}(k_s R) k_s}{2n+1}  -   \frac{j_{n-1}(k_p R) h_{n+1}(k_p R) k_p\mu }{(\lambda+2\mu)(2n+1)}  \right)\\
        & +  n R^2\rmi  \left( \frac{ -j_{n-1}(k_s R) h_{n}(k_s R) k_s^2  + j_{n-1}(k_p R) h_{n}(k_p R) k_p^2}{2n+1}   \right),
\end{split}
\]
\[
 \begin{split}
\mathfrak{c}_{2n}=& - 2(n^2-1) R \rmi  \left( \frac{j_{n+1}(k_s R) h_{n-1}(k_s R) k_s}{2n+1}  -  \frac{j_{n+1}(k_p R) h_{n-1}(k_p R) k_p\mu }{(\lambda+2\mu)(2n+1)}   \right)\\
         & -(n+1) R^2\rmi  \left( \frac{ -j_{n-1}(k_s R) h_{n}(k_s R) k_s^2   + j_{n-1}(k_p R) h_{n}(k_p R) k_p^2 }{2n+1}   \right), \\
\mathfrak{d}_{2n}=&  2(n+2)R \rmi  \left( \frac{j_{n+1}(k_s R) h_{n+1}(k_s R) k_s n}{(2n+1)}  +  \frac{j_{n+1}(k_p R) h_{n+1}(k_p R) k_p \mu (n+1)}{(\lambda+2\mu)(2n+1)}  \right)\\
         & - R^2\rmi  \left( \frac{ j_{n+1}(k_s R) h_{n}(k_s R) k_s^2 n  + j_{n+1}(k_p R) h_{n}(k_p R) k_p^2(n+1) }{2n+1}   \right).
\end{split}
\]

\end{lem}

For the static case, one has the following lemma \cite{DLL262}. 
\begin{lem}\label{lem:sinsi}
For $\bx\in  B_R$, the single layer potential has the following expression for $-1\leq m \leq 1$
\begin{equation*}
\bS_{\bbS_R}[\Tcal_1^m](\bx) =  -\frac{r}{3\mu} \Tcal_1^m,\quad
 \bS_{\bbS_R}[ \Ical_{0}^m](\bx) = -\frac{(2\lambda + 5\mu)R}{3\mu(\lambda + 2\mu)} \Ical_{0}^m.
\end{equation*}
\end{lem}

\begin{lem}\label{lem:sins}
For $\bx\in \partial B_R$, the single layer potential has the following relationship for $-1\leq m \leq 1$
\begin{equation*}
\bS_{\bbS_R}[\Tcal_1^m] = f_1 \Tcal_1^m,\quad
 \bS_{\bbS_R}[ \Ical_{0}^m] = g_{1} \Ical_{0}^m,
\end{equation*}
where
\[
f_1 = -\frac{R}{3\mu}, \quad g_1 = -\frac{(2\lambda + 5\mu)R}{3\mu(\lambda + 2\mu)}.
\]
The N-P operators satisfy for $-1\leq m \leq 1$
\[
 \bK_{\bbS_R}^{ *}[\Tcal_1^m]  = \frac{1}{2}  \Tcal_1^m, \quad \bK_{\bbS_R}^{ *}[\Ical_{0}^m]  = \frac{1}{2} \Ical_{0}^m.
\]
\end{lem}

\medskip
\subsection{The derivation of the resonant frequencies}
In this subsection, we analyze the resonant frequencies for the system \eqref{eq:mo} within the spherical configuration using two methods. One is to utilize the conclusion obtained in Theorem \ref{thm:ref}. The other is through directly solving the system \eqref{eq:mo} to obtain the resonant frequencies.

\subsubsection{The derivation of resonant frequencies from Theorem \ref{thm:ref}}
In this part, we apply the results obtained in Theorem \ref{thm:ref} to derive the resonant frequencies. To that end, we need firstly calculate the eigenvalues of the matrix $M$ defined in  Lemma \ref{lem:ma1}.
As aforementioned, the vectorial spherical harmonics $\Ical_0^m$ and $r\Tcal_1^m$, $-1\leq m\leq 1$, form a basis for the space $\mho$ given in \eqref{eq:dpsi}. Moreover, from Lemma \ref{lem:sins}, these harmonics $\Ical_0^m$ and $r\Tcal_1^m$, $-1\leq m\leq 1$, also form a basis for the space $\mho^*$. Thus the functions $\bxi_j, \bzeta_j$, $1\leq j\leq 6$ given in Lemma \ref{lem:orker} can be taken as
\[
  \bxi_j = \bzeta_j = \frac{1}{R} \Tcal_1^{j-2}, \quad 1\leq j \leq 3,
\]
and 
\[
  \bxi_j = \bzeta_j = \frac{1}{R} \Ical_0^{j-5}, \quad 4\leq j \leq 6.
\]
From Lemma \ref{lem:sinsi}, we can readily have that the matrix $M$ is a diagonal matrix, and 
\[
m_{jj}= \int_{D} -\bxi_j(\by) \cdot  \bS_{\partial D}[\bzeta_j](\by) d\by = \frac{R^2}{15\mu}, \quad 1\leq j\leq 3,
\]
as well as
\[
m_{jj}= \int_{D} -\bxi_j(\by) \cdot  \bS_{\partial D}[\bzeta_j](\by) d\by = \frac{(2\lambda + 5\mu)R^2}{9\mu(\lambda + 2\mu)} , \quad 4\leq j\leq 6.
\]
Thus the matrix $M$ has two eigenvalues
\[
 \varrho_1=\frac{R^2}{15\mu}, \quad\mbox{and}\quad  \varrho_2=\frac{(2\lambda + 5\mu)R^2}{9\mu(\lambda + 2\mu)},
\]
whose multiplicities are both $3$. Finally, the resonant frequencies from Theorem \ref{thm:ref} are 
\begin{equation}\label{eq:rethe}
\omega_1 = \frac{\epsilon^{1/2}}{R}\sqrt{\frac{15\mu}{\rho}}  \left(1 + o(1) \right) \quad \mbox{and} \quad  \omega_2 = \frac{\epsilon^{1/2}}{R}\sqrt{\frac{9 \mu(\lambda + 2\mu)}{(2\lambda + 5\mu) \rho}}  \left(1 + o(1) \right).
\end{equation}

%We first consider the case that $\bzeta = \Tcal_{1}^m$. From Lemma \ref{lem:sinsi}, the corresponding vector $\bxi = r  \Tcal_{1}^m$. Direct calculation shows that 
%\[
% \left( \bzeta ,  \bxi \right) = R^3 \int_{\mathbb{S}} \abs{ \Tcal_{1}^m}^2 ds,
%\]
%and 
%\[
% -\int_{D} \bxi(\by) \cdot  \bS_{\partial D}[\bzeta](\by) d\by = \frac{1}{3\mu}\int_{0}^R r^4 dr   \int_{\mathbb{S}} \abs{ \Tcal_{1}^m}^2 ds = \frac{R^2}{15\mu}   \int_{\mathbb{S}} \abs{ \Tcal_{1}^m}^2 ds.
%\]
%Thus the resonance frequency is  
%\[
% \omega =  \frac{\delta^{1/2}}{R}\sqrt{\frac{15\mu}{\rho}}  \left(1 + o(1) \right).
%\]
%Then, we consider the case that $\bzeta = \Ical_{0}^m$
% The corresponding vector $\bxi =  \Ical_{0}^m$. Direct calculation shows that 
%\[
% \left( \bzeta ,  \bxi \right) = R^2 \int_{\mathbb{S}} \abs{ \Ical_{0}^m}^2 ds,
%\]
%and 
%\[
% -\int_{D} \bxi(\by) \cdot  \bS_{\partial D}[\bzeta](\by) d\by =  \frac{(2\lambda + 5\mu)R}{3\mu(\lambda + 2\mu)} \int_{0}^R r^2 dr   \int_{\mathbb{S}} \abs{ \Ical_{0}^m}^2 ds = \frac{(2\lambda + 5\mu)R^4}{9\mu(\lambda + 2\mu)} \int_{\mathbb{S}} \abs{ \Ical_{0}^m}^2 ds.
%\]
%Thus the resonance frequency is  
%\[
% \omega =  \frac{\delta^{1/2}}{R}\sqrt{\frac {9\mu(\lambda + 2\mu)}{(2\lambda + 5\mu)} }  \left(1 + o(1) \right).
%\]
%

\subsubsection{The derivation of resonant frequencies by solving $(\ref{eq:conre1})$}
In this part, we derive the resonant frequencies for the system \eqref{eq:mo} through directly solving the equation \eqref{eq:conre1}. For that purpose, we first introduce the following functions:
\begin{equation*}
    \Xi_1^m = \left(
    \begin{array}{c}
     \Tcal_1^m \\
     0 \\
    \end{array}
  \right), \quad
  \Xi_2^m = \left(
    \begin{array}{c}
    0\\
     \Tcal_1^m \\
    \end{array}
  \right), \quad
  \Xi_3^m = \left(
    \begin{array}{c}
     \Ical_0^m \\
     0 \\
    \end{array}
  \right), 
\end{equation*}
and 
\begin{equation*}
 \Xi_4^m = \left(
    \begin{array}{c}
    0\\
     \Ical_0^m \\
    \end{array}
  \right), \quad
    \Xi_5^m = \left(
    \begin{array}{c}
     \Ncal_2^m \\
     0 \\
    \end{array}
  \right), \quad
  \Xi_6^m = \left(
    \begin{array}{c}
    0\\
     \Ncal_2^m \\
    \end{array}
  \right),  
\end{equation*}
where $-1\leq m\leq 1$.  Lemma \ref{lem:ker} shows that $\Xi_1^m$ and $\Xi_3^m$ are nontrivial kernels of the operator $\Acal(0, 0)$. For $\omega\ll1$, nontrivial kernels of the operator $\Acal(\omega,\delta)$ should be perturbations of $\Xi_1^m$ and $\Xi_3^m$. 
We first investigate the nontrivial kernel caused by the perturbation of the function $\Xi_1^m$. From Lemmas \ref{lem:sin} and \ref{lem:trac}, the nontrivial kernel under the perturbation of the functions $\Xi_1^m$ should be the linear combinations of $\Xi_1^m$ and $\Xi_2^m$ since the functions $\Tcal_1^m$ with $-1\leq m\leq 1$ are eigenfunctions of the operators $\bS_{\bbS_R}^{\omega}$ and $\bK_{\bbS_R}^{\omega, *}$.
%Thus one can conclude from Lemma \ref{lem:sins} that the leading order of the nontrivial solution should be $\bp \Tcal_1^m$ with $\bp=(p_1, q_2)^t$ or $\bq %\Ical_0^m$ with $\bp=(q_1, q_2)^t$. 
Denote by $M_1$ the matrix expression of the operator $\Acal(\omega, \delta)$ under the function $(\Xi_1^m, \Xi_2^m)^t$; that is 
\[
\Acal(\omega, \delta) (\Xi_1^m, \Xi_2^m) = (\Xi_1^m, \Xi_2^m) M_1.
\]
Moreover, from Lemmas \ref{lem:sin} and \ref{lem:trac}, the matrix $M_1$ has the following expression 
\[
M_1 =  \left(
    \begin{array}{cc}
       \td{b}_1 &   -{b}_1  \medskip \\
        \td{\mathfrak{b}}_1 -1 &  -{\mathfrak{b}}_1 \\
    \end{array}
  \right), 
\]
where $b_1$ and $\mathfrak{b}_1$ are given in Lemmas \ref{lem:sin} and \ref{lem:trac}, and $\tilde{b}_1$ and $\tilde{\mathfrak{b}}_1$ are given by $b_1$ and $\mathfrak{b}_1$ with $(\lambda, \mu, \rho)$ replaced by $(\title{\lambda}, \tilde{\mu}, \tilde{\rho})$. The same notations hold for other parameters $\tilde{c}_{in}, \tilde{d}_{in},  \td{\mathfrak{c}}_{in}, \td{\mathfrak{d}}_{in}$.
To ensure that there exists a nontrivial kernel for the operator $\Acal(\omega,\delta)$, the determinant of the matrix $M_1$ should vanish. From the asymptotic expansions of the functions $j_1(t)$ and $h_1(t)$ in \eqref{eq:asjt}, one has that 
\[
 \det  M_1 = 
 \left( \frac{R}{3\td{\mu}}  -   \frac{\td{\rho} R^3 \omega^2}{45\mu \tilde{\mu}} \right) \left(1 + o(1) \right) =0.
\] 
Thus the resonance frequency satisfy
\begin{equation}\label{eq:refre1}
     \omega= \frac{1}{R}\sqrt{\frac{15 {\mu}}{\td{\rho}}}  \left(1 + o(1) \right) =  \frac{\epsilon^{1/2}}{R}\sqrt{\frac{15\mu}{\rho}}  \left(1 + o(1) \right).
\end{equation}

Next, we consider the case that the resonance is caused by the perturbation of the function $\Xi_3^m$ and the situation is a bit more complex than the previous one.  From Lemmas \ref{lem:sin} and \ref{lem:trac}, the ranges of the operators $\bS_{\bbS_R}^{\omega}$ and $\bK_{\bbS_R}^{\omega, *}$ on $\Ical_0^m$  contain $\Ical_0^m$ and $\Ncal_2^m$. Thus the nontrivial kernels of the operator $\Acal(\omega,\delta)$ under the perturbation of the functions $\Xi_3^m$ should be linear combinations of $\Xi_n^m$ with $3\leq n\leq 6$. Let $\mathbf{\Xi}^m$ be defined by
\[
\mathbf{\Xi}^m=(\Xi_3^m, \Xi_4^m, \Xi_5^m, \Xi_6^m).
\]
Denote by $M_2$ the matrix expression of the operator $\Acal(\omega, \delta)$ under the function $\mathbf{\Xi}$; that is 
\[
\Acal(\omega, \delta) \mathbf{\Xi} = \mathbf{\Xi} M_2.
\]
From Lemmas \ref{lem:sin} and \ref{lem:trac}, the matrix $M_2$ has the following expression 
\[
M_2 =  \left(
    \begin{array}{cccc}
       \td{c}_{1n} &   -{c}_{1n}  &  \td{c}_{2n}  &   -{c}_{2n}  \medskip \\
       \td{\mathfrak{c}}_{1n} -1 &  -{\mathfrak{c}}_{1n}  &   \td{\mathfrak{c}}_{2n}  &  -{\mathfrak{c}}_{2n}    \medskip \\
       \td{d}_{1n} & -d_{1n} & \td{d}_{2n} &  -{d}_{2n} \medskip \\
       \td{\mathfrak{d}}_{1n}  &  -{\mathfrak{d}}_{1n}  & \td{\mathfrak{d}}_{2n} -1 &  -{\mathfrak{d}}_{2n}    \medskip
    \end{array}
  \right).
\]
To ensure the existence of a nontrivial kernel for the operator $\Acal(\omega, \delta)$, the determinant of the matrix $M_2$ should vanishes. With the help of the asymptotic expansions of the functions $j_1(t)$ and $h_1(t)$ in \eqref{eq:asjt}, one has that 
\[
\det M_2 = \frac{\delta(3\lambda + 2\mu)(\lambda + 4\mu)(2\lambda + 5\mu) R^2}{225\mu^2(\lambda + 2\mu)^3}\left( 1  -   \frac{(2\lambda + 5\mu) R^2 \rho \omega^2}{ 9 \epsilon \mu(\lambda + 2\mu)}   \right) \left(1 + o(1) \right).
\]
Thus the resonant frequency fulfill 
\begin{equation}\label{eq:refre2}
     \omega = \frac{\epsilon^{1/2}}{R}\sqrt{\frac{9 \mu(\lambda + 2\mu)}{(2\lambda + 5\mu) \rho}}  \left(1 + o(1) \right).
\end{equation}
\begin{rem}
 Resonant frequencies derived in \eqref{eq:refre1} and \eqref{eq:refre2} are in accordance with ones calculated in \eqref{eq:rethe}, which are derived via Theorem \ref{thm:ref}.
Thus, in this part, we not only provide another method to find the resonant frequencies of the system  \eqref{eq:mo} for the spherical geometry, but also validate the correctness of Theorem \ref{thm:ref}.
\end{rem}

\medskip
\subsection{The characterization of the dipolar resonances}

In this subsection, the dipolar characterization of the resonances is presented. 
It is noted that any elastic wave in a ball can be expressed by the following modes \cite{LCS5116}
%\begin{equation}\label{eq:elamo}
%\begin{split}
%\mathbf{J}_{nm1}=\nabla\left(j_n(k_p r) Y_n^m\right), \quad 0\leq n, \; -n\leq m \leq n, 
%\end{split}
%\end{equation}

\begin{equation}\label{eq:elamo}
\mathbf{J}_{nm1}=\nabla\left(j_n(k_p r) Y_n^m\right), \;\; \mathbf{J}_{nm2} = \nabla\land\left( \bx j_n(k_s r) Y_n^m\right), \;\; \mathbf{J}_{nm3} = \nabla\land\nabla\land\left( \bx j_n(k_s r) Y_n^m\right).
\end{equation}
As aforementioned, in the resonant state, the leading terms of the elastic fields inside the domain $D$ are $\bS_{\bbS_R}^{\omega}[ \Ical_{0}^m](\bx)$ and $\bS_{\bbS_R}^{\omega}[ \Tcal_1^m](\bx)$. Moreover, from Proposition \ref{pro:singlei} and by tedious calculations, one has that 
for $\bx\in B_R$, the single layer potentials $\bS_{\bbS_R}^{\omega}[ \Ical_{0}^m](\bx)$ and $\bS_{\bbS_R}^{\omega}[ \Tcal_1^m](\bx)$ have the following expression
\[
\begin{split}
 \bS_{\bbS_R}^{\omega}[ \Ical_{0}^m](\bx) =&  \frac{-\rmi  \left( 2 h_1(k_p R) + h^{\prime}(k_p R) k_p R \right)}{k_p R(\lambda + 2\mu)} \mathbf{J}_{1m1}  -\frac{\rmi  \left( 2 h_1(k_s R) + h^{\prime}(k_s R) k_s R \right)}{k_s R  \mu}   \mathbf{J}_{1m3}\\ 
 = & -\frac{(2\lambda + 5\mu)R}{3\mu(\lambda + 2\mu)} \Ical_{0}^m + \Ocal(\omega),
\end{split}
\]
and 
\[
\begin{split}
 \bS_{\bbS_R}^{\omega}[ \Tcal_1^m](\bx) &=  -\frac{\rmi k_s R^2 h_{1,s} }{\mu} \mathbf{J}_{1m2} = -\frac{r}{3\mu} \Tcal_1^m +  \Ocal(\omega).
\end{split}
\]
\begin{rem}\label{re:dipb}
It is remarked that if $n=0$, the elastic modes in \eqref{eq:elamo} denotes the monopolar wave fields, and  if $n=1$, the modes denotes the dipolar wave fields \cite{DLQ3904}.	From the last two equations, one has that the fields $\bS_{\bbS_R}^{\omega}[ \Ical_{0}^m]$ and $\bS_{\bbS_R}^{\omega}[ \Tcal_1^m]$ contain only the dipolar wave fields, $ \mathbf{J}_{1mj}$, $1\leq j \leq 3$. Moreover, in the quasi-static approximation, the leading order terms of the dipolar wave fields, $ \mathbf{J}_{1mj}$, $1\leq j \leq 3$, are $\Ical_0^m$ and $r\Tcal_1^m$, $-1\leq m\leq 1$. As aforementioned, the vectorial spherical harmonics $\Ical_0^m$ and $r\Tcal_1^m$, $-1\leq m\leq 1$, forms a basis for the space $\mho$ given in \eqref{eq:dpsi}. Thus the basis functions $\{\bxi_j\}^6_{j=1}$ of $\mho$ can be treated as dipolar wave fields in the static regime. 
\end{rem}

\begin{rem}	
The conclusions discussed above can be extended to the two-dimensional case via applying spectral properties of the single-layer potential operators and N-P operator \cite{LLZ2555}.
\end{rem}

\section*{Acknowledgment}
%The work of H. Liu is supported by the Hong Kong RGC General Research Funds (projects 11311122, 11300821 and 12301420),  the NSFC/RGC Joint Research Fund (project N\_CityU101/21), and the ANR/RGC Joint Research Grant, A\_CityU203/19. 
The work of Hongjie Li was supported by the NSF/RGC Joint Research Fund (project N\_CityU101/21).
The work of J. Zou is supported by the Hong Kong RGC General Research Funds (projects 14306921 and 14306719).

 \bibliographystyle{abbrv}
\bibliography{ela_res_pap}{}

\end{document}